\documentclass[leqno,12pt]{amsart}
\usepackage{amssymb} 
\theoremstyle{plain} 
\newtheorem{theorem}{\indent\sc Theorem}[section] 
\newtheorem{lemma}[theorem]{\indent\sc Lemma}
\newtheorem{corollary}[theorem]{\indent\sc Corollary}
\newtheorem{proposition}[theorem]{\indent\sc Proposition}

\newtheorem{question}[theorem]{\indent\sc Question}
\theoremstyle{definition} 
\newtheorem{definition}[theorem]{\indent\sc Definition}

\newtheorem{example}[theorem]{\indent\sc Example}

\begin{document}

\title{Studies on the Chazy equations \\}
\author{Yusuke Sasano }

\renewcommand{\thefootnote}{\fnsymbol{footnote}}
\footnote[0]{2000\textit{ Mathematics Subjet Classification}.
34M55; 34M45; 58F05; 32S65.}

\keywords{ 
Birational symmetry, Chazy equations, Painlev\'e equations.}
\maketitle

\begin{abstract}
In this paper, we study the Chazy III,IX and X equations. For the Chazy III equation, by making the birational transformations
the Chazy III equation is transformed into a third-order ordinary differential equation of rational type. For this equation, we find its meromorphic solutions, whose free parameters are essentially two. We also show that the system associated with this equation admits new special solutions solved by $tanh(t)$. For the Chazy IX equation, we transform the Chazy IX equation to a system of the first-order ordinary differential equations by birational transformations. For this system, we give two new birational B{\"a}cklund transformations. We also give the holomorphy condition of this system. Thanks to this holomorphy condition, we obtain a new partial differential system in two variables involving the Chazy IX equation, This system satisfies the compatibility condition, and admits a travelling wave solution. For the Chazy X equation, we transform the Chazy X equation to a system of the first-order ordinary differential equations by birational transformations. For this system, we give two birational B{\"a}cklund transformations. One of them is new. We also give the holomorphy condition of this system. Thanks to this holomorphy condition, we can recover this system.
\end{abstract}

\section{Introduction}

In 1910, Chazy studied Painlev\'e type equation with third order (see \cite{6,Cos1}) explicitly given by
\begin{equation}\label{S1}
\frac{d^3 u}{dt^3}=2u\frac{d^2u}{dt^2}-3\left(\frac{du}{dt}\right)^2.
\end{equation}
Here $u$ denotes unknown complex variable. It is known that this equation fails some Painlev\'e test \cite{1,6,Cos1}. Nevertheless, Chazy gave this special attention \cite{1,2,7}.

This equation has a solution \cite{1}
\begin{equation}
u(t)=4\frac{d}{dt}\rm{log\rm}{\theta_1}'(0,t),
\end{equation}
and special solutions
\begin{equation}
  \left\{
  \begin{aligned}
   u_1(t) &=\frac{a}{(t-t_0)^2}-\frac{6}{t-t_0} \quad (a,t_0 \in {\Bbb C}),\\
   u_2(t) &=-\frac{1}{t-t_0} \quad (t_0 \in {\Bbb C}).\\
   \end{aligned}
  \right. 
\end{equation}

\begin{figure}[t]
\unitlength 0.1in
\begin{picture}(57.90,24.60)(11.90,-24.80)
%
\special{pn 8}%
\special{pa 2000 1200}%
\special{pa 3200 220}%
\special{fp}%
\special{pa 2020 1190}%
\special{pa 4170 1190}%
\special{fp}%
%
\special{pn 8}%
\special{pa 2010 1210}%
\special{pa 3220 2390}%
\special{fp}%
%
\special{pn 8}%
\special{pa 3200 2370}%
\special{pa 4160 1190}%
\special{dt 0.045}%
\special{pa 4160 1190}%
\special{pa 4160 1191}%
\special{dt 0.045}%
%
\special{pn 8}%
\special{pa 3200 230}%
\special{pa 4160 1180}%
\special{dt 0.045}%
\special{pa 4160 1180}%
\special{pa 4160 1180}%
\special{dt 0.045}%
%
\special{pn 8}%
\special{pa 3210 240}%
\special{pa 3210 2350}%
\special{dt 0.045}%
\special{pa 3210 2350}%
\special{pa 3210 2349}%
\special{dt 0.045}%
%
\special{pn 8}%
\special{pa 5870 230}%
\special{pa 4670 1190}%
\special{fp}%
%
\special{pn 8}%
\special{pa 5880 240}%
\special{pa 6980 1180}%
\special{fp}%
%
\special{pn 8}%
\special{pa 4690 1180}%
\special{pa 6960 1180}%
\special{dt 0.045}%
\special{pa 6960 1180}%
\special{pa 6959 1180}%
\special{dt 0.045}%
%
\special{pn 8}%
\special{pa 4680 1180}%
\special{pa 5870 2410}%
\special{dt 0.045}%
\special{pa 5870 2410}%
\special{pa 5870 2410}%
\special{dt 0.045}%
\special{pa 5870 2410}%
\special{pa 6960 1180}%
\special{dt 0.045}%
\special{pa 6960 1180}%
\special{pa 6960 1180}%
\special{dt 0.045}%
%
\special{pn 8}%
\special{pa 5870 260}%
\special{pa 5870 2410}%
\special{dt 0.045}%
\special{pa 5870 2410}%
\special{pa 5870 2409}%
\special{dt 0.045}%
%
\special{pn 20}%
\special{sh 0.600}%
\special{ar 2020 1190 20 20  0.0000000 6.2831853}%
%
\special{pn 20}%
\special{sh 0.600}%
\special{ar 2880 920 20 20  0.0000000 6.2831853}%
%
\special{pn 20}%
\special{sh 0.600}%
\special{ar 5700 890 20 20  0.0000000 6.2831853}%
%
\special{pn 20}%
\special{sh 0.600}%
\special{ar 5880 230 20 20  0.0000000 6.2831853}%
%
\special{pn 20}%
\special{sh 0.600}%
\special{ar 6360 640 20 20  0.0000000 6.2831853}%
%
\special{pn 20}%
\special{pa 3200 2380}%
\special{pa 2900 2090}%
\special{fp}%
\special{sh 1}%
\special{pa 2900 2090}%
\special{pa 2934 2151}%
\special{pa 2938 2127}%
\special{pa 2962 2122}%
\special{pa 2900 2090}%
\special{fp}%
%
\special{pn 20}%
\special{pa 5870 2400}%
\special{pa 5870 2060}%
\special{fp}%
\special{sh 1}%
\special{pa 5870 2060}%
\special{pa 5850 2127}%
\special{pa 5870 2113}%
\special{pa 5890 2127}%
\special{pa 5870 2060}%
\special{fp}%
\put(25.1000,-24.0000){\makebox(0,0)[lb]{$u:=X$}}%
\put(56.5000,-20.6000){\makebox(0,0)[lb]{$v:=Y$}}%
\put(17.5000,-11.2000){\makebox(0,0)[lb]{$u_1(t)$}}%
\put(26.7000,-8.7000){\makebox(0,0)[lb]{$u_2(t)$}}%
\put(55.4000,-1.9000){\makebox(0,0)[lb]{$v_1(t)$}}%
\put(62.1000,-5.6000){\makebox(0,0)[lb]{$v_2(t)$}}%
\put(54.6000,-8.3000){\makebox(0,0)[lb]{$v_3(t)$}}%
\put(11.9000,-16.8000){\makebox(0,0)[lb]{${\Bbb P}^3$}}%
\put(23.9000,-26.5000){\makebox(0,0)[lb]{Classical equation $\eqref{S1}$}}%
\put(52.0000,-26.5000){\makebox(0,0)[lb]{New equation $\eqref{SSSS1}$}}%
\end{picture}%
\label{fig:Chazy4}
\caption{Each figure denotes 3-dimensional projective space $(X,Y,Z) \in {\Bbb C}^3 \subset {\Bbb P}^3$. The solutions $u_1,u_2$ and $v_3$ are rational solutions which pass through an accessible singular point (see Section 3) in the boundary divisor ${\mathcal H}$, respectively. The solutions $v_1,v_2$ are new meromorphic solutions which pass through an accessible singular point in the boundary divisor ${\mathcal H}$, respectively. The equation \eqref{S1} can be written by the variable $u$ and new equation \eqref{SSSS1} can be written by the variable $v:=-\frac{\frac{du}{dt}}{u}$. For the equation \eqref{S1} we can not find the solutions with initial conditions in the part surrounding the dotted line. The solutions $u_2(t),v_3(t)$ are common in each equation.}
\end{figure}
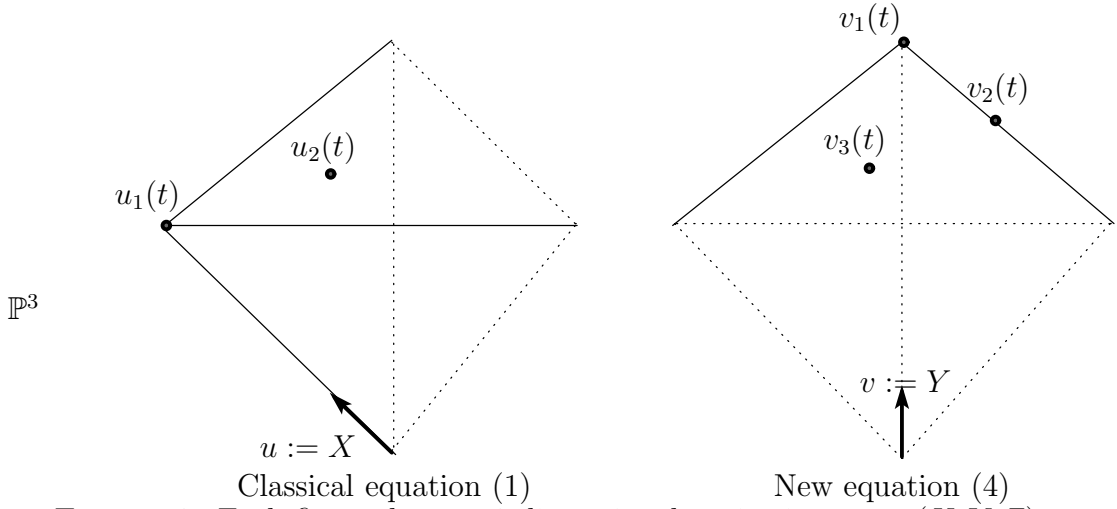

In this paper, by making a change of variables
\begin{equation*}
X=\frac{u}{6}, \quad Y=-\frac{\frac{du}{dt}}{u}, \quad Z=-\frac{\frac{d^2u}{dt^2}}{\frac{du}{dt}}+\frac{u}{3},
\end{equation*}
the system \eqref{S1} is transformed as follows:
\begin{equation*}
\frac{dX}{dt}=-XY, \quad \frac{dY}{dt}=(2X+Y-Z)Y, \quad \frac{dZ}{dt}=Z^2+8XZ-20XY-20X^2.
\end{equation*}
By elimination of $X,Z$ and setting $v:=Y$, we obtain
\begin{equation}\label{SSSS1}
\frac{d^3v}{dt^3}=-3\left(v^2-\frac{dv}{dt}\right)\left(\frac{dv}{dt}\right)+\frac{3}{2}v^4-\frac{\left(v^3-2\frac{d^2v}{dt^2}\right) \left(5v^3+2\frac{d^2v}{dt^2}\right)}{2\left(v^2+2\frac{dv}{dt}\right)}.
\end{equation}
We find its meromorphic solutions given by
\begin{equation}
  \left\{
  \begin{aligned}
   v_1(t) &=\frac{-1}{t-t_0}+a_1+a_2(t-t_0)+\frac{a_1(a_1^2+a_2)}{2}(t-t_0)^2+\cdots,\\
   v_2(t) &=\frac{-2}{t-t_0}+a_1(t-t_0)^2+\frac{2}{21}a_1^2(t-t_0)^5+\cdots,\\
   v_3(t) &=\frac{1}{t-t_0},\\
   \end{aligned}
  \right. 
\end{equation}
where $a_1,a_2$ are free parameters and $t_0$ is an initial position. The solutions $v_1(t),v_2(t)$ are new (see Figure 1). The solutions $u_2(t),v_3(t)$ are common in each system.

We also show that this system admits new special solutions (see Section 5):
\begin{equation}
(X,Y,Z)=(c_1,0,-6c_1tanh(6(c_1 t-c_1c_2))-4c_1) \quad (c_1,c_2 \in {\Bbb C}).
\end{equation}

We remark that the classical Darboux-Halphen system (see \cite{25})
\begin{equation}
  \left\{
  \begin{aligned}
   \frac{dx}{dt} &=yz-x(y+z),\\
   \frac{dy}{dt} &=xz-y(x+z),\\
   \frac{dz}{dt} &=xy-z(x+y)
   \end{aligned}
  \right. 
\end{equation}
is equivalent to the equation \eqref{S1} if one sets
\begin{equation}
u=-2(x+y+z), \quad \frac{du}{dt}=2(xz+yz+xy), \quad \frac{d^2u}{dt^2}=-12xyz.
\end{equation}
This system is invariant under the transformation:
\begin{equation*}
\pi:(x,y,z) \rightarrow (y,z,x) \quad (\pi^3=1).
\end{equation*}

It is well-known that this system has the following rational solutions:
\begin{equation*}
  \left\{
  \begin{aligned}
   x_1(t) &=\frac{1}{(t-t_0)},\\
   y_1(t) &=\frac{1}{(t-t_0)},\\
   z_1(t) &=\frac{1}{(t-t_0)} \quad (t_0 \in {\Bbb C}),
   \end{aligned}
  \right. 
\end{equation*}
and
\begin{equation*}
  \left\{
  \begin{aligned}
   x_2(t) &=\frac{1}{(t-t_0)},\\
   y_2(t) &=\frac{1}{(t-t_0)},\\
   z_2(t) &=\frac{a}{(t-t_0)^2}+\frac{b}{(t-t_0)} \quad (a,b,t_0 \in {\Bbb C}).
   \end{aligned}
  \right. 
\end{equation*}

Reviews of Chazy's work on III and results of further research can be found in \cite{1,2,7}. Chazy-III, or an equivalent system of first-order equations, appears in several physics contexts, for example, self-dual Yang-Mills equations \cite{3}. Clarkson and Olver \cite{7} obtained III among the group-invariant reductions of the partial differential equation,
\begin{equation*}
w_{xxx}=w_{y}w_{xx}-w_{x}w_{yy},
\end{equation*}
which has applications in boundary-layer theory. These authors and C. M. Cosgrove also gave a theory of higher-order equations having properties similar to III.

We also study the Chazy IX equation:
\begin{equation}
\frac{d^3 u}{dt^3}=54u^4+72u^2\frac{du}{dt}+12\left(\frac{du}{dt} \right)^2+\delta \quad (\delta \in {\Bbb C}),
\end{equation}
where $u$ denotes unknown complex variable and $\delta$ is its constant parameter.

In 2000, C. M. Cosgrove constructed the general solution of the Chazy IX equation in terms of hyperelliptic functions (see \cite{Cos1}).

For this equation, we give two birational B{\"a}cklund transformations:
\begin{align}
\begin{split}
&g_0(u;\delta) \rightarrow \\
&\left(\frac{(\sqrt{5}-3)\{108u^4+18(5+\sqrt{5})u^2 u'+6(3+\sqrt{5})(u')^2+3(\sqrt{5}-1)uu''+2\delta}{6(\sqrt{5}-1)\{3(\sqrt{5}-1)uu'+u''\}};\frac{7+3\sqrt{5}}{2} \delta \right),\\
&g_1(u;\delta) \rightarrow \\
&\left(\frac{(-\sqrt{5}-3)\{108u^4+18(5-\sqrt{5})u^2 u'+6(3-\sqrt{5})(u')^2+3(-\sqrt{5}-1)uu''+2\delta}{6(-\sqrt{5}-1)\{3(-\sqrt{5}-1)uu'+u''\}};\frac{7-3\sqrt{5}}{2} \delta \right),
\end{split}
\end{align}
where $u'=\frac{du}{dt},\ u''=\frac{d^2u}{dt^2}$. These B{\"a}cklund transformations are new.

We show that the birational transformation
\begin{equation}
  \left\{
  \begin{aligned}
   x &=u,\\
   y &=\frac{du}{dt}+\frac{3}{2}(\sqrt{5}-1)u^2,\\
   z &=\frac{d^2u}{dt^2}+3(\sqrt{5}-1)u \frac{du}{dt}
   \end{aligned}
  \right. 
\end{equation}
takes the Chazy IX equation to the system of the first-order ordinary differential equations:
\begin{equation}
  \left\{
  \begin{aligned}
   \frac{dx}{dt} &=-\frac{3}{2}(\sqrt{5}-1)x^2+y,\\
   \frac{dy}{dt} &=z,\\
   \frac{dz}{dt} &=3(\sqrt{5}+3)y^2+3(\sqrt{5}-1)xz+\delta.
   \end{aligned}
  \right. 
\end{equation}
We make this system in the polynomial class from the viewpoint of accessible singularity and local index.

For this system, we give the holomorphy condition of this system. Thanks to this holomorphy condition, we obtain a new partial differential system in two variables $(t,s)$ involving the Chazy IX equation:
\begin{equation}\label{eq:AA1}
  \left\{
  \begin{aligned}
   dx =&\left\{-\frac{3}{2}(\sqrt{5}-1)x^2+y\right\}dt\\
   &+\left\{-12(\sqrt{5}+2)x^2y+(\sqrt{5}+1)\left(3xz-2y^2-\frac{2}{3}\delta \right) \right\}ds,\\
   dy =&z dt\\
   &+\{-6(2\sqrt{5}-5)x^2z+12\sqrt{5}xy^2-(\sqrt{5}+1)yz+(3\sqrt{5}-5)\delta x \}ds,\\
   dz =&\{3(\sqrt{5}+3)y^2+3(\sqrt{5}-1)xz+\delta\}dt\\
   &+\{48xyz-24y^3-(\sqrt{5}+1)z^2+2(\sqrt{5}-3)\delta y \}ds.
   \end{aligned}
  \right. 
\end{equation}
We remark that when $s=0$, we can obtain the Chazy IX system.

We will show that this system satisfies the compatibility condition, and admits a travelling wave solution.

Moreover,we find the relation between the system \eqref{eq:AA1} and soliton equations (see Section 11). In this paper, we can make the birational transformations between the system \eqref{eq:AA1} and the partial differential system:
\begin{equation}\label{eq:AA4}
  \left\{
  \begin{aligned}
   \frac{\partial^3 u}{\partial t^3} =&54u^4+72u^2 \frac{\partial u}{\partial t}+12\left( \frac{\partial u}{\partial t} \right)^2+\delta,\\
   \frac{\partial u}{\partial S} =&54u^4+18u^2 \frac{\partial u}{\partial t}+3\left( \frac{\partial u}{\partial t} \right)^2-\frac{9}{2}u\frac{\partial^2 u}{\partial t^2}+\delta,
   \end{aligned}
  \right. 
\end{equation}
where $u:=x$. The first equation in \eqref{eq:AA4} coincides with Chazy IX equation.

Adding each system in \eqref{eq:AA4}, we can obtain
\begin{equation}\label{eq:AA5}
   \frac{\partial^3 u}{\partial t^3} =54u^2\frac{\partial u}{\partial t}+9\left( \frac{\partial u}{\partial t} \right)^2+\frac{9}{2}u \frac{\partial^2 u}{\partial t^2}+\frac{\partial u}{\partial S}.
\end{equation}

\begin{question}
It is still an open question whether the second equation in \eqref{eq:AA4} and the system \eqref{eq:AA5} coincide with which of the soliton equations.
\end{question}

We also study the Chazy X equation (see \cite{Cos1}):
\begin{align}
\begin{split}
X.a:\frac{d^3 u}{dt^3}=&6u^2 \frac{du}{dt}+\frac{3}{11}(9+7\sqrt{3})\left(\frac{du}{dt}+u^2 \right)^2\\
&-\frac{1}{22}(4-3\sqrt{3})\alpha \frac{du}{dt}+\frac{1}{44}(3-5\sqrt{3})\alpha u^2-\frac{1}{352}(9+7\sqrt{3}){\alpha}^2,\\
X.b:\frac{d^3 u}{dt^3}=&6u^2 \frac{du}{dt}+\frac{3}{11}(9-7\sqrt{3})\left(\frac{du}{dt}+u^2 \right)^2\\
&-\frac{1}{22}(4+3\sqrt{3})\alpha \frac{du}{dt}+\frac{1}{44}(3+5\sqrt{3})\alpha u^2-\frac{1}{352}(9-7\sqrt{3}){\alpha}^2,
\end{split}
\end{align}
where $u$ denotes unknown complex variable and $\alpha$ is its constant parameter.

The full version of Class X with $\alpha \not=0$ was found by C. M. Cosgrove (see \cite{Cos1}), the parameter $\alpha$ being overlooked in Chazy's original work. Chazy looked for recessive terms in his Class X but missed them because of some minor calculation error. Bureau (see \cite{Bureau}), having derived the 13 reduced Chazy classes from first principles, did not attempt to attach recessive terms to this class. Exton (see \cite{Exton}), who incorrectly concluded that Chazy Classes IX and X were unstable equations, also apparantly did not look for resessive terms.

The splitting of Chazy Class X deserves a comment. Because its coefficients contain the irrationality $\sqrt{3}$ in an essential way, changing the sign of $\sqrt{3}$ yields a distinct equation. The two versions, which we denote Chazy-X.a and Chazy-X.b, are related by a B{\"a}cklund transformation, which was also found by C. M. Cosgrove (see \cite{Cos1}).

In 2000, C. M. Cosgrove constructed the general solution of the Chazy X equation in terms of hyperelliptic functions (see \cite{Cos1}).

For these equations, we give two birational B{\"a}cklund transformations from Chazy-X.a to Chazy-X.b:
\begin{align}
\begin{split}
&g_0(u;\alpha) \rightarrow \\
&(\frac{3+\sqrt{3}}{1056\{(3+\sqrt{3})uu'-u''\}} \{288(6+\sqrt{3})u^2 u'-176(3+\sqrt{3})uu''+96(9+7\sqrt{3})u^4\\
&+64(-3+5\sqrt{3})(u')^2-8(-3+5\sqrt{3})\alpha u^2+16(-4+3\sqrt{3})\alpha u'-(9+7\sqrt{3})\alpha^2 \};\alpha ),\\
&g_1(u;\alpha) \rightarrow \\
&\left(\frac{(4+\sqrt{3})\{8(-4+\sqrt{3})u^3-8(-4+\sqrt{3})u u'-8(-4+\sqrt{3})u''-(-7+5\sqrt{3})\alpha u \}}{13\{8u^2+8u'+(\sqrt{3}-1)\alpha\}};(-2+\sqrt{3})\alpha \right),
\end{split}
\end{align}
where $u'=\frac{du}{dt},\ u''=\frac{d^2u}{dt^2}$. The transformation $g_1$ is new. Thanks to these B{\"a}cklund transformations, in this paper we discuss the case of Chazy-X.a equation.

We show that the birational transformation
\begin{equation}
  \left\{
  \begin{aligned}
   x &=u,\\
   y &=\frac{du}{dt}-\frac{3+\sqrt{3}}{2}u^2,\\
   z &=\frac{d^2u}{dt^2}-(3+\sqrt{3})u \frac{du}{dt}
   \end{aligned}
  \right. 
\end{equation}
takes the Chazy X.a equation to the system of the first-order ordinary differential equations:
\begin{equation}
  \left\{
  \begin{aligned}
   \frac{dx}{dt} &=\frac{3+\sqrt{3}}{2}x^2+y,\\
   \frac{dy}{dt} &=z,\\
   \frac{dz}{dt} &=\frac{2}{11}(-3+5\sqrt{3})y^2-(3+\sqrt{3})xz+\frac{1}{22}(-4+3\sqrt{3})\alpha y-\frac{1}{352}(9+7\sqrt{3})\alpha^2.
   \end{aligned}
  \right. 
\end{equation}
We make this system in the polynomial class from the viewpoint of accessible singularity and local index.

For this system, we give the holomorphy condition of this system. Thanks to this holomorphy condition, we can recover the Chazy X.a equation.

\section{The Chazy polynomial class}
In \cite{6}, Chazy attempted the complete classification of all third-order differential equations of the form:
\begin{equation}
\frac{d^3u}{dt^3}=F\left(t,u,\frac{du}{dt},\frac{d^2u}{dt^2}\right),
\end{equation}
where $F$ is a polynomial in $u,\frac{du}{dt}$ and $\frac{d^2u}{dt^2}$ and locally analytic in $t$, having the Painlev\'e property.

Chazy proved that Painlev\'e-type differential equations of the third-order in the polynomial class must take the form:
\begin{align}
\begin{split}
\frac{d^3u}{dt^3}=&Qu\frac{d^2u}{dt^2}+R\left(\frac{du}{dt}\right)^2+Su^2\frac{du}{dt}+Tu^4+A(t)\frac{d^2u}{dt^2}+B(t)u\frac{du}{dt}\\
&+C(t)\frac{du}{dt}+D(t)u^3+E(t)u^2+F(t)u+G(t),
\end{split}
\end{align}
where, after a suitable normalization, $Q,R,S$ and $T$ are certain rational or algebraic numbers, and the remaining coefficients are locally analytic functions of the complex variable $t$ to be determined.

The canonical reduced equations defining each of the Chazy classes is as follows:
\begin{align}
\begin{split}
I:&\frac{d^3u}{dt^3}=-6\left(\frac{du}{dt}\right)^2;\\
II:&\frac{d^3u}{dt^3}=-2u\frac{d^2u}{dt^2}-2\left(\frac{du}{dt}\right)^2;\\
III:&\frac{d^3u}{dt^3}=2u\frac{d^2u}{dt^2}-3\left(\frac{du}{dt}\right)^2;\\
IV:&\frac{d^3u}{dt^3}=-3u\frac{d^2u}{dt^2}-3\left(\frac{du}{dt}\right)^2-3u^2\frac{du}{dt};\\
V:&\frac{d^3u}{dt^3}=-2u\frac{d^2u}{dt^2}-4\left(\frac{du}{dt}\right)^2-2u^2\frac{du}{dt};\\
VI:&\frac{d^3u}{dt^3}=-u\frac{d^2u}{dt^2}-5\left(\frac{du}{dt}\right)^2-u^2\frac{du}{dt};\\
VII:&\frac{d^3u}{dt^3}=-u\frac{d^2u}{dt^2}-2\left(\frac{du}{dt}\right)^2+2u^2\frac{du}{dt};\\
VIII:&\frac{d^3u}{dt^3}=6u^2\frac{du}{dt};\\
IX:&\frac{d^3u}{dt^3}=12\left(\frac{du}{dt}\right)^2+72u^2\frac{du}{dt}+54u^4;\\
X.a:&\frac{d^3u}{dt^3}=6u^2\frac{du}{dt}+\frac{3}{11}(9+7\sqrt{3})\left(\frac{du}{dt}+u^2\right)^2;\\
X.b:&\frac{d^3u}{dt^3}=6u^2\frac{du}{dt}+\frac{3}{11}(9-7\sqrt{3})\left(\frac{du}{dt}+u^2\right)^2;
\end{split}
\end{align}
\begin{align*}
XI:&\frac{d^3u}{dt^3}=-2u\frac{d^2u}{dt^2}-2\left(\frac{du}{dt}\right)^2+\frac{24}{N^2-1}\left(\frac{du}{dt}+u^2\right)^2;\\
XII:&\frac{d^3u}{dt^3}=2u\frac{d^2u}{dt^2}-3\left(\frac{du}{dt}\right)^2-\frac{4}{N^2-36}\left(6\frac{du}{dt}-u^2\right)^2;\\
XIII:&\frac{d^3u}{dt^3}=12u\frac{du}{dt}.
\end{align*}
In Chazy-XI, $N$ is a positive integer not equal to 1 or a multiple of 6. In Chazy-XII, $N$ is a positive integer not equal to 1 or 6.

It is well-known that the KdV equation belongs to Chazy Class XIII, the modified KdV equation belongs to Class VIII, and the potential KdV equation (as well as the soliton equation $u_{xxt}=-6u_xu_t+\alpha$) belongs to Class I.

\section{Accessible singularity and local index}
Let us review the notion of {\it accessible singularity}. Let $B$ be a connected open domain in $\Bbb C$ and $\pi : {\mathcal W} \longrightarrow B$ a smooth proper holomorphic map. We assume that ${\mathcal H} \subset {\mathcal W}$ is a normal crossing divisor which is flat over $B$. Let us consider a rational vector field $\tilde v$ on $\mathcal W$ satisfying the condition
\begin{equation*}
\tilde v \in H^0({\mathcal W},\Theta_{\mathcal W}(-\log{\mathcal H})({\mathcal H})).
\end{equation*}
Fixing $t_0 \in B$ and $P \in {\mathcal W}_{t_0}$, we can take a local coordinate system $(x_1,\ldots ,x_n)$ of ${\mathcal W}_{t_0}$ centered at $P$ such that ${\mathcal H}_{\rm smooth \rm}$ can be defined by the local equation $x_1=0$.
Since $\tilde v \in H^0({\mathcal W},\Theta_{\mathcal W}(-\log{\mathcal H})({\mathcal H}))$, we can write down the vector field $\tilde v$ near $P=(0,\ldots ,0,t_0)$ as follows:
\begin{equation*}
\tilde v= \frac{\partial}{\partial t}+g_1 
\frac{\partial}{\partial x_1}+\frac{g_2}{x_1} 
\frac{\partial}{\partial x_2}+\cdots+\frac{g_n}{x_1} 
\frac{\partial}{\partial x_n}.
\end{equation*}
This vector field defines the following system of differential equations
\begin{equation}\label{39}
\frac{dx_1}{dt}=g_1(x_1,\ldots,x_n,t),\ \frac{dx_2}{dt}=\frac{g_2(x_1,\ldots,x_n,t)}{x_1},\cdots, \frac{dx_n}{dt}=\frac{g_n(x_1,\ldots,x_n,t)}{x_1}.
\end{equation}
Here $g_i(x_1,\ldots,x_n,t), \ i=1,2,\dots ,n,$ are holomorphic functions defined near $P=(0,\dots ,0,t_0).$

\begin{definition}\label{Def1}
With the above notation, assume that the rational vector field $\tilde v$ on $\mathcal W$ satisfies the condition
$$
(A) \quad \tilde v \in H^0({\mathcal W},\Theta_{\mathcal W}(-\log{\mathcal H})({\mathcal H})).
$$
We say that $\tilde v$ has an {\it accessible singularity} at $P=(0,\dots ,0,t_0)$ if
$$
x_1=0 \ {\rm and \rm} \ g_i(0,\ldots,0,t_0)=0 \ {\rm for \rm} \ {\rm every \rm} \ i, \ 2 \leq i \leq n.
$$
\end{definition}

If $P \in {\mathcal H}_{{\rm smooth \rm}}$ is not an accessible singularity, all solutions of the ordinary differential equation passing through $P$ are vertical solutions, that is, the solutions are contained in the fiber ${\mathcal W}_{t_0}$ over $t=t_0$. If $P \in {\mathcal H}_{\rm smooth \rm}$ is an accessible singularity, there may be a solution of \eqref{39} which passes through $P$ and goes into the interior ${\mathcal W}-{\mathcal H}$ of ${\mathcal W}$.

Here we review the notion of {\it local index}. Let $v$ be an algebraic vector field with an accessible singular point $\overrightarrow{p}=(0,\ldots,0)$ and $(x_1,\ldots,x_n)$ be a coordinate system in a neighborhood centered at $\overrightarrow{p}$. Assume that the system associated with $v$ near $\overrightarrow{p}$ can be written as
\begin{align}\label{b}
\begin{split}
\frac{d}{dt}\begin{pmatrix}
             x_1 \\
             x_2 \\
             \vdots\\
             x_{n-1} \\
             x_n
             \end{pmatrix}=\frac{1}{x_1}\left\{\begin{bmatrix}
             a_{11} & 0 & 0 & \hdots & 0 \\
             a_{21} & a_{22} & 0 &  \hdots & 0 \\
             \vdots & \vdots & \ddots & 0 & 0 \\
             a_{(n-1)1} & a_{(n-1)2} & \hdots & a_{(n-1)(n-1)} & 0 \\
             a_{n1} & a_{n2} & \hdots & a_{n(n-1)} & a_{nn}
             \end{bmatrix}\begin{pmatrix}
             x_1 \\
             x_2 \\
             \vdots\\
             x_{n-1} \\
             x_n
             \end{pmatrix}+\begin{pmatrix}
             x_1h_1(x_1,\ldots,x_n,t) \\
             h_2(x_1,\ldots,x_n,t) \\
             \vdots\\
             h_{n-1}(x_1,\ldots,x_n,t) \\
             h_n(x_1,\ldots,x_n,t)
             \end{pmatrix}\right\},\\
              (h_i \in {\Bbb C}(t)[x_1,\ldots,x_n], \ a_{ij} \in {\Bbb C}(t))
             \end{split}
             \end{align}
where $h_1$ is a polynomial which vanishes at $\overrightarrow{p}$ and $h_i$, $i=2,3,\ldots,n$ are polynomials of order at least 2 in $x_1,x_2,\ldots,x_n$, We call ordered set of the eigenvalues $(a_{11},a_{22},\cdots,a_{nn})$ {\it local index} at $\overrightarrow{p}$.

We are interested in the case with local index
\begin{equation}\label{integer}
(1,a_{22}/a_{11},\ldots,a_{nn}/a_{11}) \in {\Bbb Z}^{n}.
\end{equation}
These properties suggest the possibilities that $a_1$ is the residue of the formal Laurent series:
\begin{equation}
y_1(t)=\frac{a_{11}}{(t-t_0)}+b_1+b_2(t-t_0)+\cdots+b_n(t-t_0)^{n-1}+\cdots \quad (b_i \in {\Bbb C}),
\end{equation}
and the ratio $(1,a_{22}/a_{11},\ldots,a_{nn}/a_{11})$ is resonance data of the formal Laurent series of each $y_i(t) \ (i=2,\ldots,n)$, where $(y_1,\ldots,y_n)$ is original coordinate system satisfying $(x_1,\ldots,x_n)=(f_1(y_1,\ldots,y_n),\ldots,f_n(y_1,\ldots,y_n)), \ f_i(y_1,\ldots,y_n) \in {\Bbb C}(t)(y_1,\ldots,y_n)$.

If each component of $(1,a_{22}/a_{11},\ldots,a_{nn}/a_{11})$ has the same sign, we may resolve the accessible singularity by blowing-up finitely many times. However, when different signs appear, we may need to both blow up and blow down.

The $\alpha$-test,
\begin{equation}\label{poiuy}
t=t_0+\alpha T, \quad x_i=\alpha X_i, \quad \alpha \rightarrow 0,
\end{equation}
yields the following reduced system:
\begin{align}\label{ppppppp}
\begin{split}
\frac{d}{dT}\begin{pmatrix}
             X_1 \\
             X_2 \\
             \vdots\\
             X_{n-1} \\
             X_n
             \end{pmatrix}=\frac{1}{X_1}\begin{bmatrix}
             a_{11}(t_0) & 0 & 0 & \hdots & 0 \\
             a_{21}(t_0) & a_{22}(t_0) & 0 &  \hdots & 0 \\
             \vdots & \vdots & \ddots & 0 & 0 \\
             a_{(n-1)1}(t_0) & a_{(n-1)2}(t_0) & \hdots & a_{(n-1)(n-1)}(t_0) & 0 \\
             a_{n1}(t_0) & a_{n2}(t_0) & \hdots & a_{n(n-1)}(t_0) & a_{nn}(t_0)
             \end{bmatrix}\begin{pmatrix}
             X_1 \\
             X_2 \\
             \vdots\\
             X_{n-1} \\
             X_n
             \end{pmatrix},
             \end{split}
             \end{align}
where $a_{ij}(t_0) \in {\Bbb C}$. Fixing $t=t_0$, this system is the system of the first order ordinary differential equation with constant coefficient. Let us solve this system. At first, we solve the first equation:
\begin{equation}
X_1(T)=a_{11}(t_0)T+C_1 \quad (C_1 \in {\Bbb C}).
\end{equation}
Substituting this into the second equation in \eqref{ppppppp}, we can obtain the first order linear ordinary differential equation:
\begin{equation}
\frac{dX_2}{dT}=\frac{a_{22}(t_0) X_2}{a_{11}(t_0)T+C_1}+a_{21}(t_0).
\end{equation}
By variation of constant, in the case of $a_{11}(t_0) \not= a_{22}(t_0)$ we can solve explicitly:
\begin{equation}
X_2(T)=C_2(a_{11}(t_0)T+C_1)^{\frac{a_{22}(t_0)}{a_{11}(t_0)}}+\frac{a_{21}(t_0)(a_{11}(t_0)T+C_1)}{a_{11}(t_0)-a_{22}(t_0)} \quad (C_2 \in {\Bbb C}).
\end{equation}
This solution is a single-valued solution if and only if
$$
\frac{a_{22}(t_0)}{a_{11}(t_0)} \in {\Bbb Z}.
$$
In the case of $a_{11}(t_0)=a_{22}(t_0)$ we can solve explicitly:
\begin{equation}
X_2(T)=C_2(a_{11}(t_0)T+C_1)+\frac{a_{21}(t_0)(a_{11}(t_0)T+C_1){\rm Log}(a_{11}(t_0)T+C_1)}{a_{11}(t_0)} \quad (C_2 \in {\Bbb C}).
\end{equation}
This solution is a single-valued solution if and only if
$$
a_{21}(t_0)=0.
$$
Of course, $\frac{a_{22}(t_0)}{a_{11}(t_0)}=1 \in {\Bbb Z}$.
In the same way, we can obtain the solutions for each variables $(X_3,\ldots,X_n)$. The conditions $\frac{a_{jj}(t)}{a_{11}(t)} \in {\Bbb Z}, \ (j=2,3,\ldots,n)$ are necessary condition in order to have the Painlev\'e property.

\begin{example}
For an example of the condition \eqref{integer}, let us consider
\begin{equation}\label{ex1}
\frac{d^3 u}{dt^3}=u\frac{d^2u}{dt^2}-2\left(\frac{du}{dt}\right)^2+6u^2\frac{du}{dt}.
\end{equation}
Here $u$ denotes unknown complex variable. We will show that this equation is the integrable non-Painlev\'e equation.

\begin{proposition}
The birational transformation
\begin{equation}
  \left\{
  \begin{aligned}
   x &=u,\\
   y &=-\frac{\frac{du}{dt}}{u}+u,\\
   z &=-\frac{\frac{d^2u}{dt^2}}{\frac{du}{dt}}+2u
   \end{aligned}
  \right. 
\end{equation}
takes the equation \eqref{ex1} to the system
\begin{equation}\label{system;ex1}
  \left\{
  \begin{aligned}
   \frac{dx}{dt} &=x^2-xy,\\
   \frac{dy}{dt} &=y^2-xy+xz-yz,\\
   \frac{dz}{dt} &=z^2-3xz-4xy.
   \end{aligned}
  \right. 
\end{equation}
Here $x,y,z$ denote unknown complex variables.
\end{proposition}

We will show that this system violates the condition \eqref{integer}.

Let us take the coordinate system $(p,q,r)$ centered at the point $(p,q,r)=(0,0,0)$:
$$
p=\frac{1}{x}, \quad q=\frac{y}{x}, \quad r=\frac{z}{x}.
$$
The system \eqref{system;ex1} is rewritten as follows:
\begin{align*}
\frac{d}{dt}\begin{pmatrix}
             p \\
             q \\
             r 
             \end{pmatrix}&=\frac{1}{p}\left\{\begin{pmatrix}
             -1 & 0 & 0 \\
             0 & -2 & 1 \\
             0 & -4 & -4
             \end{pmatrix}\begin{pmatrix}
             p \\
             q \\
             r 
             \end{pmatrix}+\cdots\right\}
             \end{align*}
satisfying \eqref{b}. To the above system, we make the linear transformation
\begin{equation*}
\begin{pmatrix}
             X \\
             Y \\
             Z 
             \end{pmatrix}=\begin{pmatrix}
             1 & 0 & 0 \\
             0 & \frac{2\sqrt{-3}}{3} & \frac{3+\sqrt{-3}}{6} \\
             0 & -\frac{2\sqrt{-3}}{3} & \frac{3-\sqrt{-3}}{6}
             \end{pmatrix}\begin{pmatrix}
             p \\
             q \\
             r 
             \end{pmatrix}
\end{equation*}
to arrive at
\begin{equation*}
\frac{d}{dt}\begin{pmatrix}
             X \\
             Y \\
             Z 
             \end{pmatrix}=\frac{1}{X}\left\{\begin{pmatrix}
             -1 & 0 & 0 \\
             0 & -3+\sqrt{-3} & 0 \\
             0 & 0 & -3-\sqrt{-3}
             \end{pmatrix}\begin{pmatrix}
             X \\
             Y \\
             Z 
             \end{pmatrix}+\cdots\right\}.
             \end{equation*}
             In this case, the ratio $(\frac{-3+\sqrt{-3}}{-1},\frac{-3-\sqrt{-3}}{-1})=(3-\sqrt{-3},3+\sqrt{-3})$ is not in ${\Bbb Z}^2$.

\end{example}

\begin{example}
For an application of the condition \eqref{integer}, let us consider
\begin{equation}\label{ex2}
\frac{d^3 u}{dt^3}=2u\frac{d^2u}{dt^2}-3\left(\frac{du}{dt}\right)^2+a\left(6\frac{du}{dt}-u^2\right)^2 \quad (a \in {\Bbb C}).
\end{equation}
Here $u$ denotes unknown complex variable.

\begin{proposition}
The birational transformation
\begin{equation}
  \left\{
  \begin{aligned}
   x &=\frac{u^2}{\frac{du}{dt}},\\
   y &=\frac{\frac{du}{dt}}{u},\\
   z &=\frac{\frac{d^2u}{dt^2}}{\frac{du}{dt}}
   \end{aligned}
  \right. 
\end{equation}
takes the equation \eqref{ex2} to the system
\begin{equation}\label{system;ex2}
  \left\{
  \begin{aligned}
   \frac{dx}{dt} &=2xy-xz,\\
   \frac{dy}{dt} &=-y^2+yz,\\
   \frac{dz}{dt} &=-3xy^2+36axy^2-12ax^2y^2+ax^3y^2+2xyz-z^2.
   \end{aligned}
  \right. 
\end{equation}
Here $x,y,z$ denote unknown complex variables.
\end{proposition}

We will consider when this system satisfies the condition \eqref{integer}.

Let us take the coordinate system $(p,q,r)$ centered at the point $(p,q,r)=(0,0,0)$:
$$
p=x-\frac{3a-\sqrt{9a^2-a}}{a}, \quad q=\frac{y}{z}-\frac{1}{2}, \quad r=\frac{1}{z}.
$$
The system \eqref{system;ex2} is rewritten as follows:
\begin{align*}
\frac{d}{dt}\begin{pmatrix}
             p \\
             q \\
             r 
             \end{pmatrix}&=\frac{1}{p}\left\{\begin{pmatrix}
             0 & -\frac{2(-3a+\sqrt{9a^2-a})}{a} & 0 \\
             \frac{1-9a-3\sqrt{9a^2-a}}{4} & -\frac{a+2\sqrt{9a^2-a}}{2a} & 0 \\
             0 & 0 & -\frac{1}{2}
             \end{pmatrix}\begin{pmatrix}
             p \\
             q \\
             r 
             \end{pmatrix}+\cdots\right\}
             \end{align*}
satisfying \eqref{b}. To the above system, we make the linear transformation
\begin{equation*}
\begin{pmatrix}
             X \\
             Y \\
             Z 
             \end{pmatrix}=\begin{pmatrix}
             -\frac{\sqrt{a(9a-1)}}{4-42a+14\sqrt{a(9a-1)}} & \frac{a}{a-2\sqrt{a(9a-1)}} & 0 \\
             \frac{\sqrt{a(9a-1)}}{4-42a+14\sqrt{a(9a-1)}} & \frac{2-18a}{2-18a+\sqrt{a(9a-1)}} & 0 \\
             0 & 0 & 1
             \end{pmatrix}\begin{pmatrix}
             p \\
             q \\
             r 
             \end{pmatrix}
\end{equation*}
to arrive at
\begin{equation*}
\frac{d}{dt}\begin{pmatrix}
             X \\
             Y \\
             Z 
             \end{pmatrix}=\frac{1}{Z}\left\{\begin{pmatrix}
             -\frac{1}{2} & 0 & 0 \\
             0 & -\frac{\sqrt{a(9a-1)}}{a} & 0 \\
             0 & 0 & -\frac{1}{2}
             \end{pmatrix}\begin{pmatrix}
             X \\
             Y \\
             Z 
             \end{pmatrix}+\cdots\right\}.
             \end{equation*}
             In this case, the ratio $\left(\frac{-\frac{1}{2}}{-\frac{1}{2}},\frac{-\frac{\sqrt{a(9a-1)}}{a}}{-\frac{1}{2}} \right)$ is in ${\Bbb Z}^2$ if and only if
\begin{equation}
4a(9a-1)=N^2 a^2 \quad (N \in {\Bbb Z}).
\end{equation}
This equation can be solved as follows:
\begin{equation}
a=-\frac{4}{N^2-36} \quad (N \in {\Bbb N}).
\end{equation}
This case coincides with Chazy-XII equation.

\end{example}

In the next section, in order to consider the phase spaces for each system, let us take the compactification $[z_0:z_1:z_2:z_3] \in {\Bbb P}^3$ of $(x,y,z) \in {\Bbb C}^3$ with the natural embedding
$$
(x,y,z)=(z_1/z_0,z_2/z_0,z_3/z_0).
$$
Moreover, we denote the boundary divisor in ${\Bbb P}^3$ by $ {\mathcal H}$. Extend the regular vector field on ${\Bbb C}^3$ to a rational vector field $\tilde v$ on ${\Bbb P}^3$. It is easy to see that ${\Bbb P}^3$ is covered by four copies of ${\Bbb C}^3${\rm : \rm}
\begin{align*}
&U_0={\Bbb C}^3 \ni (x,y,z),\\
&U_j={\Bbb C}^3 \ni (X_j,Y_j,Z_i) \ (j=1,2,3),
\end{align*}
via the following rational transformations
\begin{align}\label{P2cover}
\begin{split}
& X_1=1/x, \quad Y_1=y/x, \quad Z_1=z/x,\\
& X_2=x/y, \quad Y_2=1/y, \quad Z_2=z/y,\\
& X_3=x/z, \quad Y_3=y/z, \quad Z_3=1/z.
\end{split}
\end{align}

\section{Chazy-III equation}
Chazy-III equation is given by
\begin{equation*}
\frac{d^3 u}{dt^3}=2u\frac{d^2u}{dt^2}-3\left(\frac{du}{dt}\right)^2.
\end{equation*}
Here $u$ denotes unknown complex variable.

\begin{proposition}
The birational transformation
\begin{equation}
  \left\{
  \begin{aligned}
   x &=\frac{u}{6},\\
   y &=-\frac{\frac{du}{dt}}{u}+\frac{u}{6},\\
   z &=-\frac{\frac{d^2u}{dt^2}}{\frac{du}{dt}}+\frac{u}{3}
   \end{aligned}
  \right. 
\end{equation}
takes the equation \eqref{S1} to the system
\begin{equation}\label{1}
  \left\{
  \begin{aligned}
   \frac{dx}{dt} &=x^2-xy,\\
   \frac{dy}{dt} &=y^2-xy+xz-yz,\\
   \frac{dz}{dt} &=z^2+8xz-20xy.
   \end{aligned}
  \right. 
\end{equation}
Here $x,y,z$ denote unknown complex variables.
\end{proposition}

Let us do the Painlev\'e test. To find the leading order behaviour of a singularity at $t=t_0$ one sets
\begin{equation*}
  \left\{
  \begin{aligned}
   x & \propto \frac{a_0}{(t-t_0)},\\
   y & \propto \frac{b_0}{(t-t_0)},\\
   z & \propto \frac{c_0}{(t-t_0)},
   \end{aligned}
  \right. 
\end{equation*}
from which it is easily deduced that
\begin{enumerate}
\item $(a_0,b_0,c_0)=(-1,0,0),$
\item $(a_0,b_0,c_0)=(0,-1,0),$
\item $(a_0,b_0,c_0)=(0,0,-1),$
\item $(a_0,b_0,c_0)=(0,-2,-1).$
\end{enumerate}

{\bf Case $(a_0,b_0,c_0)=(-1,0,0)$}

In this case, we find
\begin{equation*}
  \left\{
  \begin{aligned}
   x(t) &=\frac{-1}{(t-t_0)} \quad (t_0 \in {\Bbb C}),\\
   y(t) &=0,\\
   z(t) &=0.
   \end{aligned}
  \right. 
\end{equation*}

{\bf Case $(a_0,b_0,c_0)=(0,-1,0)$}

In this case, we find
\begin{equation*}
  \left\{
  \begin{aligned}
   x(t) &=a_2(t-t_0)-\frac{a_2c_1}{2}(t-t_0)^2+\cdots,\\
   y(t) &=\frac{-1}{(t-t_0)}+\frac{c_1}{2}+\frac{(28a_2+c_1^2)}{4}(t-t_0)+\frac{(8a_2+c_1^2)c_1}{8}(t-t_0)^2+\cdots,\\
   z(t) &=c_1+(20a_2+c_1^2)(t-t_0)+(14a_2+c_1^2)c_1(t-t_0)^2+\cdots,
   \end{aligned}
  \right. 
\end{equation*}
where $(a_2,c_1)$ are free parameters.

{\bf Case $(a_0,b_0,c_0)=(0,0,-1)$}

In this case, we find
\begin{equation*}
  \left\{
  \begin{aligned}
   x(t) &=a_1-\frac{a_1b_2}{2}(t-t_0)^2+\cdots,\\
   y(t) &=a_1+b_2(t-t_0)+\frac{11a_1b_2}{2}(t-t_0)^2+\cdots,\\
   z(t) &=\frac{-1}{(t-t_0)}-4a_1-12a_1^2(t-t_0)-4a_1b_2(t-t_0)^2+\cdots,
   \end{aligned}
  \right. 
\end{equation*}
where $(a_1,b_2)$ are free parameters.

{\bf Case $(a_0,b_0,c_0)=(0,-2,-1)$}

In this case, we find
\begin{equation*}
  \left\{
  \begin{aligned}
   x(t) &=a_3(t-t_0)^2-\frac{4a_3^2}{5}(t-t_0)^5+\cdots,\\
   y(t) &=\frac{-2}{(t-t_0)}+\frac{17a_3}{5}(t-t_0)^2-\frac{44a_3^2}{175}(t-t_0)^5+\cdots,\\
   z(t) &=\frac{-1}{(t-t_0)}+8a_3(t-t_0)^2+\frac{172a_3^2}{35}(t-t_0)^5+\cdots,
   \end{aligned}
  \right. 
\end{equation*}
where $a_3$ is a free parameter.

\begin{proposition}
The system \eqref{1} becomes again a system in the polynomial class in each coordinate system{\rm:\rm}
\begin{center}
${U_j}={\Bbb C}^3 \ni \{(x_j,y_j,z_j)\},  \ \ j=0,1,\ldots,4$
\end{center}
via the following birational transformations{\rm:\rm}
\begin{align}
\begin{split}
0) \ &x_0=x, \quad y_0=y, \quad z_0=z,\\
1) \ &x_1=-\frac{x-y}{2x}, \quad y_1=x, \quad z_1=\frac{(x-y)(x+3y-2z)}{4x},\\
2) \ &x_2=xy, \quad y_2=\frac{1}{y}, \quad z_2=z,\\
3) \ &x_3=x, \quad y_3=(y-x)z, \quad z_3=\frac{1}{z},\\
4) \ &x_4=-(x-y)xz, \quad y_4=-\frac{1}{(x-y)z}, \quad z_4=\frac{1}{z}.
\end{split}
\end{align}
\end{proposition}

\begin{corollary}
The system \eqref{1} is equivalent to the following systems\rm{:\rm}
\begin{equation}
  \left\{
  \begin{aligned}
   \frac{dx_1}{dt} &=x_1^2y_1-z_1,\\
   \frac{dy_1}{dt} &=-2x_1y_1^2,\\
   \frac{dz_1}{dt} &=2y_1(6x_1^3y_1+5x_1z_1+6z_1),
   \end{aligned}
  \right. 
\end{equation}
and
\begin{equation}
  \left\{
  \begin{aligned}
   \frac{dx_2}{dt} &=(x_2y_2^2-1)x_2z_2,\\
   \frac{dy_2}{dt} &=-(x_2y_2^2-1)(y_2z_2-1),\\
   \frac{dz_2}{dt} &=z_2^2+8x_2y_2z_2-20x_2,
   \end{aligned}
  \right. 
\end{equation}
and
\begin{equation}
  \left\{
  \begin{aligned}
   \frac{dx_3}{dt} &=-x_3y_3z_3,\\
   \frac{dy_3}{dt} &=y_3(-20x_3y_3z_3^2-20x_3^2z_3+y_3z_3+10x_3),\\
   \frac{dz_3}{dt} &=-1-8x_3z_3+20x_3^2z_3^2+20x_3y_3z_3^3,
   \end{aligned}
  \right. 
\end{equation}
and
\begin{equation}
  \left\{
  \begin{aligned}
   \frac{dx_4}{dt} &=-10x_4^2(2x_4y_4^2z_4+2z_4^2-y_4),\\
   \frac{dy_4}{dt} &=20x_4^2y_4^3z_4+20x_4y_4z_4^2-10x_4y_4^2-z_4,\\
   \frac{dz_4}{dt} &=-1-8x_4y_4z_4+20x_4^2y_4^2z_4^2+20x_4z_4^3.
   \end{aligned}
  \right. 
\end{equation}
\end{corollary}

The following Lemma shows that the rational vector field $\tilde v$ associated with the system \eqref{1} has six accessible singular points on the boundary divisor ${\mathcal H} \subset {\Bbb P}^3$.
\begin{lemma}\label{9.1}
The rational vector field $\tilde v$ has six accessible singular points{\rm : \rm}
\begin{equation}
  \left\{
  \begin{aligned}
   P_1 &=\{(X_1,Y_1,Z_1)|X_1=Y_1=Z_1=0\},\\
   P_2 &=\{(X_2,Y_2,Z_2)|X_2=Y_2=Z_2=0\},\\
   P_3 &=\{(X_3,Y_3,Z_3)|X_3=Y_3=Z_3=0\},\\
   P_4 &=\{(X_1,Y_1,Z_1)|X_1=0, \ Y_1=1, \ Z_1=2\},\\
   P_5 &=\{(X_1,Y_1,Z_1)|X_1=0, \ Y_1=1, \ Z_1=-10\},\\
   P_6 &=\{(X_2,Y_2,Z_2)|X_2=Y_2=0, \ Z_2=\frac{1}{2}\},\\
   \end{aligned}
  \right. 
\end{equation}
where $P_4$ is multiple point of order $2$.
\end{lemma}
This lemma can be proven by a direct calculation. \qed

Next let us calculate its local index at each point.
\begin{center}
\begin{tabular}{|c|c|c|} \hline 
Singular point & Type of local index   \\ \hline 
$P_1$ & $(-1,3,2)$  \\ \hline 
$P_2$ & $(2,1,1)$  \\ \hline 
$P_3$ & $(1,2,1)$  \\ \hline 
$P_5$ & $(0,12,-12)$  \\ \hline 
$P_6$ & $(3,1,-2)$  \\ \hline 
\end{tabular}
\end{center}

\begin{example}
Let us take the coordinate system $(p,q,r)$ centered at the point $P_1$:
$$
p=X_1=\frac{1}{x}, \quad q=Y_1=\frac{y}{x}, \quad r=Z_1=\frac{z}{x}.
$$
The system \eqref{1} is rewritten as follows:
\begin{align*}
\frac{d}{dt}\begin{pmatrix}
             p \\
             q \\
             r 
             \end{pmatrix}&=\frac{1}{p}\left\{\begin{pmatrix}
             -1 & 0 & 0 \\
             0 & -2 & 1 \\
             0 & -20 & 7 
             \end{pmatrix}\begin{pmatrix}
             p \\
             q \\
             r 
             \end{pmatrix}+\cdots \right\}
             \end{align*}
satisfying \eqref{b}. To the above system, we make the linear transformation
\begin{equation*}
\begin{pmatrix}
             X \\
             Y \\
             Z 
             \end{pmatrix}=\begin{pmatrix}
             1 & 0 & 0 \\
             0 & -4 & 1 \\
             0 & 5 & -1 
             \end{pmatrix}\begin{pmatrix}
             p \\
             q \\
             r 
             \end{pmatrix}
\end{equation*}
to arrive at
\begin{equation*}
\frac{d}{dt}\begin{pmatrix}
             X \\
             Y \\
             Z 
             \end{pmatrix}=\frac{1}{X}\left\{\begin{pmatrix}
             -1 & 0 & 0 \\
             0 & 3 & 0 \\
             0 & 0 & 2 
             \end{pmatrix}\begin{pmatrix}
             X \\
             Y \\
             Z 
             \end{pmatrix}+\cdots \right\}.
             \end{equation*}
             In this case, the local index is $(-1,3,2)$. This suggests the possibilities that $-1$ is the residue of the formal Laurent series:
\begin{equation}
x(t)=\frac{-1}{(t-t_0)}+a_1+a_2(t-t_0)+\dots+a_n(t-t_0)^{n-1}+\cdots \quad (a_i \in {\Bbb C}),
\end{equation}
and the ratio $(\frac{3}{-1},\frac{2}{-1})=(-3,-2)$ is resonance data of the formal Laurent series of $(y(t),z(t))$, respectively. We see that the formal Laurent series which passes through $P_1$ have no free parameters. There is only one solution which passes through $P_1$ explicitly given by
             \begin{equation}
             x(t)=-\frac{1}{(t-t_0)}, \ y(t)=0, \ z(t)=0 \quad (t_0 \in {\Bbb C}).
             \end{equation}
             This is a rational solution.
\end{example}

\begin{example}
Let us take the coordinate system $(p,q,r)$ centered at the point $P_5$:
$$
p=X_1=\frac{1}{x}, \quad q=Y_1-1=\frac{y}{x}-1, \quad r=Z_1+10=\frac{z}{x}+10.
$$
The system \eqref{1} is rewritten as follows:
\begin{align*}
\frac{d}{dt}\begin{pmatrix}
             p \\
             q \\
             r 
             \end{pmatrix}&=\frac{1}{p}\left\{\begin{pmatrix}
             0 & 0 & 0 \\
             0 & 12 & 0 \\
             0 & -30 & -12 
             \end{pmatrix}\begin{pmatrix}
             p \\
             q \\
             r 
             \end{pmatrix}+\cdots\right\}
             \end{align*}
satisfying \eqref{b}. To the above system, we make the linear transformation
\begin{equation*}
\begin{pmatrix}
             X \\
             Y \\
             Z 
             \end{pmatrix}=\begin{pmatrix}
             1 & 0 & 0 \\
             0 & -\frac{1}{4} & 0 \\
             0 & \frac{5}{4} & 1 
             \end{pmatrix}\begin{pmatrix}
             p \\
             q \\
             r 
             \end{pmatrix}
\end{equation*}
to arrive at
\begin{equation*}
\frac{d}{dt}\begin{pmatrix}
             X \\
             Y \\
             Z 
             \end{pmatrix}=\frac{1}{X}\left\{\begin{pmatrix}
             0 & 0 & 0 \\
             0 & 12 & 0 \\
             0 & 0 & -12 
             \end{pmatrix}\begin{pmatrix}
             X \\
             Y \\
             Z 
             \end{pmatrix}+\cdots\right\}.
             \end{equation*}
             In this case, the local index is $(0,12,-12)$. We see that the residue of the formal Laurent series:
\begin{equation}
x(t)=\frac{a_0}{(t-t_0)}+a_1+a_2(t-t_0)+\dots+a_n(t-t_0)^{n-1}+\cdots \quad (a_i \in {\Bbb C})
\end{equation}
             is equal to $a_0=0$. By a direct calculation, we see that there are no solutions which pass through $P_5$.
\end{example}

\begin{example}
Let us take the coordinate system $(p,q,r)$ centered at the point $P_3$:
$$
p=X_3=\frac{x}{z}, \quad q=Y_3=\frac{y}{z}, \quad r=Z_3=\frac{1}{z}.
$$
The system \eqref{1} is rewritten as follows:
\begin{align*}
\frac{d}{dt}\begin{pmatrix}
             p \\
             q \\
             r 
             \end{pmatrix}&=\frac{1}{r}\left\{\begin{pmatrix}
             -1 & 0 & 0 \\
             1 & -2 & 0 \\
             0 & 0 & -1 
             \end{pmatrix}\begin{pmatrix}
             p \\
             q \\
             r 
             \end{pmatrix}+\cdots\right\}
             \end{align*}
satisfying \eqref{b}. To the above system, we make the linear transformation
\begin{equation*}
\begin{pmatrix}
             X \\
             Y \\
             Z 
             \end{pmatrix}=\begin{pmatrix}
             1 & 0 & 0 \\
             1 & 1 & 0 \\
             0 & 0 & 1 
             \end{pmatrix}\begin{pmatrix}
             p \\
             q \\
             r 
             \end{pmatrix}
\end{equation*}
to arrive at
\begin{equation*}
\frac{d}{dt}\begin{pmatrix}
             X \\
             Y \\
             Z 
             \end{pmatrix}=\frac{1}{Z}\left\{\begin{pmatrix}
             -1 & 0 & 0 \\
             0 & -2 & 0 \\
             0 & 0 & -1 
             \end{pmatrix}\begin{pmatrix}
             X \\
             Y \\
             Z 
             \end{pmatrix}+\cdots\right\}.
             \end{equation*}
             In this case, the local index is $(-1,-2,-1)$. This suggests the possibilities that $-1$ is the residue of the formal Laurent series:
\begin{equation}
z(t)=\frac{-1}{(t-t_0)}+c_1+c_2(t-t_0)+\cdots+c_n(t-t_0)^{n-1}+\cdots \quad (c_i \in {\Bbb C}),
\end{equation}
and the ratio $(\frac{-1}{-1},\frac{-2}{-1})=(1,2)$ is resonance data of the formal Laurent series of $(x(t),y(t))$, respectively. There exist meromorphic solutions with three free parameters which passes through $P_3$.
\end{example}

\begin{example}
Let us take the coordinate system $(p,q,r)$ centered at the point $P_6$:
$$
p=X_2=\frac{x}{y}, \quad q=Y_2=\frac{1}{y}, \quad r=Z_2-\frac{1}{2}=\frac{z}{y}-\frac{1}{2}.
$$
The system \eqref{1} is rewritten as follows:
\begin{align*}
\frac{d}{dt}\begin{pmatrix}
             p \\
             q \\
             r 
             \end{pmatrix}&=\frac{1}{q}\left\{\begin{pmatrix}
             -\frac{3}{2} & 0 & 0 \\
             0 & -\frac{1}{2} & 0 \\
             -\frac{63}{4} & 0 & 1 
             \end{pmatrix}\begin{pmatrix}
             p \\
             q \\
             r 
             \end{pmatrix}+\cdots\right\}
             \end{align*}
satisfying \eqref{b}. To the above system, we make the linear transformation
\begin{equation*}
\begin{pmatrix}
             X \\
             Y \\
             Z 
             \end{pmatrix}=\begin{pmatrix}
             \frac{63}{10} & 0 & 0 \\
             0 & 1 & 0 \\
             -\frac{63}{10} & 0 & 1 
             \end{pmatrix}\begin{pmatrix}
             p \\
             q \\
             r 
             \end{pmatrix}
\end{equation*}
to arrive at
\begin{equation*}
\frac{d}{dt}\begin{pmatrix}
             X \\
             Y \\
             Z 
             \end{pmatrix}=\frac{1}{Y}\left\{\begin{pmatrix}
             -\frac{3}{2} & 0 & 0 \\
             0 & -\frac{1}{2} & 0 \\
             0 & 0 & 1 
             \end{pmatrix}\begin{pmatrix}
             X \\
             Y \\
             Z 
             \end{pmatrix}+\cdots\right\}.
             \end{equation*}
             In this case, the local index is $(-\frac{3}{2},-\frac{1}{2},1)$. This suggests the possibilities that $-\frac{1}{2}$ is the residue of the formal Laurent series:
\begin{equation}
y(t)=\frac{-\frac{1}{2}}{(t-t_0)}+b_1+b_2(t-t_0)+\dots+b_n(t-t_0)^{n-1}+\cdots \quad (b_i \in {\Bbb C}),
\end{equation}
and the ratio $(\frac{-\frac{3}{2}}{-\frac{1}{2}},\frac{1}{-\frac{1}{2}})=(3,-2)$ is resonance data of the formal Laurent series of $(x(t),z(t))$, respectively. There exist meromorphic solutions with two free parameters which passes through $P_6$.
\end{example}

\section{Particular solutions of the system \eqref{1}}
We see that the system \eqref{1} admits a particular solution $x=0$. Moreover $(y,z)$ satisfy
\begin{equation}
  \left\{
  \begin{aligned}
   \frac{dy}{dt} &=y^2-yz,\\
   \frac{dz}{dt} &=z^2.
   \end{aligned}
  \right. 
\end{equation}
The equation $\frac{dz}{dt}=z^2$ can be solved as follows:
\begin{equation}
z(t)=-\frac{1}{t+c_1} \quad (c_1 \in {\Bbb C}).
\end{equation}
By substituting this solution to the equation $\frac{dy}{dt}=y^2-yz$, we obtain
\begin{equation}
\frac{dy}{dt}=y^2+\frac{y}{t+c_1}.
\end{equation}
This equation can be solved by
\begin{equation}
y(t)=-\frac{2(t+c_1)}{t^2+2c_1 t-2c_2} \quad (c_1,c_2 \in {\Bbb C}).
\end{equation}

We also see that by making a change of variables
\begin{equation*}
X:=x, \quad Y:=y-x, \quad Z:=z
\end{equation*}
the system \eqref{1} is transformed as follows:
\begin{equation}\label{ChazyIIIsys}
  \left\{
  \begin{aligned}
   \frac{dX}{dt} &=-XY,\\
   \frac{dY}{dt} &=(2X+Y-Z)Y,\\
   \frac{dZ}{dt} &=Z^2+8XZ-20XY-20X^2.
   \end{aligned}
  \right. 
\end{equation}
By elimination of $X,Z$ and setting $v:=Y$, we obtain
\begin{equation*}
\frac{d^3v}{dt^3}=-3\left(v^2-\frac{dv}{dt}\right)\left(\frac{dv}{dt}\right)+\frac{3}{2}v^4-\frac{\left(v^3-2\frac{d^2v}{dt^2}\right) \left(5v^3+2\frac{d^2v}{dt^2}\right)}{2\left(v^2+2\frac{dv}{dt}\right)}.
\end{equation*}

This system admits a particular solution $Y=0$. Moreover $(X,Z)$ satisfy
\begin{equation}
  \left\{
  \begin{aligned}
   \frac{dX}{dt} &=0,\\
   \frac{dZ}{dt} &=Z^2+8XZ-20X^2.
   \end{aligned}
  \right. 
\end{equation}
By substituting $X=c_1 \ (c_1 \in {\Bbb C})$ to the equation $\frac{dZ}{dt}=Z^2+8XZ-20X^2$, we obtain
\begin{equation}
\frac{dZ}{dt}=Z^2+8c_1Z-20c_1^2.
\end{equation}
This system can be solved by
\begin{equation}
Z(t)=-6c_1tanh(6(c_1 t-c_1c_2))-4c_1 \quad (c_2 \in {\Bbb C}).
\end{equation}

\section{The Chazy-IX equation}
Chazy-IX equation is given by
\begin{equation}\label{eq;ChazyIX}
\frac{d^3 u}{dt^3}=54u^4+72u^2\frac{du}{dt}+12\left(\frac{du}{dt} \right)^2+\delta \quad (\delta \in {\Bbb C}).
\end{equation}
In this paper, at first we transform the equation \eqref{eq;ChazyIX} to a system of differential equations by birational transformations. For this system, we give two new B{\"a}cklund transformations. We also give the holomorphy condition of this system. Thanks to this condition, we obtain a new partial differential system in two variables $(t,s)$ involving this system, This system satisfies the compatibility condition, and admits a travelling wave solution.

\begin{theorem}\label{th:7.1}
The birational transformation $\varphi_0$
\begin{equation}
  \left\{
  \begin{aligned}
   x &=u,\\
   y &=\frac{du}{dt}+\frac{3}{2}(\sqrt{5}-1)u^2,\\
   z &=\frac{d^2u}{dt^2}+3(\sqrt{5}-1)u \frac{du}{dt}
   \end{aligned}
  \right. 
\end{equation}
takes the equation \eqref{eq;ChazyIX} to the system
\begin{equation}\label{system;ChazyIX}
  \left\{
  \begin{aligned}
   \frac{dx}{dt} &=-\frac{3}{2}(\sqrt{5}-1)x^2+y,\\
   \frac{dy}{dt} &=z,\\
   \frac{dz}{dt} &=3(\sqrt{5}+3)y^2+3(\sqrt{5}-1)xz+\delta.
   \end{aligned}
  \right. 
\end{equation}
\end{theorem}
Before we will prove Theorem \ref{th:7.1}, we review the case of the second Painlev\'e equation.

\section{The case of the second Painlev\'e system}

In this section, we review the case of the second Painlev\'e system:
\begin{equation}\label{PII}
\frac{d^2u}{dt^2}=2u^3+tu+\alpha \quad (\alpha \in {\Bbb C}).
\end{equation}

Let us make its polynomial Hamiltonian from the viewpoint of accessible singularity and local index.

{\bf Step 0:} We make a change of variables.
\begin{equation}
x=u, \quad y=\frac{du}{dt}.
\end{equation}

{\bf Step 1:} We make a change of variables.
\begin{equation}
x_1=\frac{1}{x}, \quad y_1=\frac{y}{x^2}.
\end{equation}
In this coordinate system, we see that this system has two accessible singular points:
\begin{equation}
(x_1,y_1)=\left\{(0,1),(0,-1) \right\}.
\end{equation}

Around the point $(x_1,y_1)=(0,1)$, we can rewrite the system as follows.

{\bf Step 2:} We make a change of variables.
\begin{equation}
x_2=x_1, \quad y_2=y_1-1.
\end{equation}
In this coordinate system, we can rewrite the system satisfying the condition \eqref{b}:
\begin{align*}
\frac{d}{dt}\begin{pmatrix}
             x_2 \\
             y_2 
             \end{pmatrix}&=\frac{1}{x_2}\left\{\begin{pmatrix}
             -1 & 0   \\
             0 & -4 
             \end{pmatrix}\begin{pmatrix}
             x_2 \\
             y_2 
             \end{pmatrix}+\cdots\right\},
             \end{align*}
and we can obtain the local index $(-1,-4)$ at the point $\{(x_2,y_2)=(0,0)\}$. The ratio of the local index at the point $\{(x_2,y_2)=(0,0)\}$ is a positive integer.

We aim to obtain the local index $(-1,-2)$ by successive blowing-up procedures.

{\bf Step 3:} We blow up at the point $\{(x_2,y_2)=(0,0)\}$.
\begin{equation}
x_3=x_2, \quad y_3=\frac{y_2}{x_2}.
\end{equation}

{\bf Step 4:} We blow up at the point $\{(x_3,y_3)=(0,0)\}$.
\begin{equation}
x_4=x_3, \quad y_4=\frac{y_3}{x_3}.
\end{equation}
In this coordinate system, we see that this system has the following accessible singular point:
\begin{equation}
(x_4,y_4)=(0,t/2).
\end{equation}

{\bf Step 5:} We make a change of variables.
\begin{equation}
x_5=x_4, \quad y_5=y_4-t/2.
\end{equation}
In this coordinate system, we can rewrite the system as follows:
\begin{align*}
\frac{d}{dt}\begin{pmatrix}
             x_5 \\
             y_5 
             \end{pmatrix}&=\frac{1}{x_5}\left\{\begin{pmatrix}
             -1 & 0  \\
             \alpha-1/2 & -2
             \end{pmatrix}\begin{pmatrix}
             x_5 \\
             y_5 
             \end{pmatrix}+\cdots\right\},
             \end{align*}
and we can obtain the local index $(-1,-2)$. Here, the relation between $(x_5,y_5)$ and $(x,y)$ is given by
\begin{equation*}
  \left\{
  \begin{aligned}
   x_5 &=\frac{1}{x},\\
   y_5 &=y-x^2-\frac{t}{2}.
   \end{aligned}
  \right. 
\end{equation*}

Finally, we can choose canonical variables $(q,p)$.

{\bf Step 9:} We make a change of variables.
\begin{equation}
q=\frac{1}{x_5}, \quad p=y_5,
\end{equation}
and we can obtain the system
\begin{equation*}
  \left\{
  \begin{aligned}
   \frac{dq}{dt} &=q^2+p+\frac{t}{2},\\
   \frac{dp}{dt} &=-2qp+\alpha-\frac{1}{2}
   \end{aligned}
  \right. 
\end{equation*}
with the polynomial Hamiltonian $H_{II}$:
\begin{equation}
H_{II}=q^2 p+\frac{1}{2}p^2+\frac{t}{2}p-\left(\alpha-\frac{1}{2} \right)q.
\end{equation}
We remark that we can discuss the case of the accessible singular point $(x_1,y_1)=(0,-1)$ in the same way as in the case of $(x_1,y_1)=(0,1)$.

\section{Proof of Theorem \ref{th:7.1}}
By the same way of the second Painlev\'e system, we can prove Theorem \ref{th:7.1}.

{\bf Proof.} At first, we rewrite the equation \eqref{eq;ChazyIX} to the system of the first-order ordinary differential equations.

{\bf Step 0:} We make a change of variables.
\begin{equation}
x=u, \quad y=\frac{du}{dt}, \quad z=\frac{d^2u}{dt^2}.
\end{equation}

{\bf Step 1:} We make a change of variables.
\begin{equation}
x_1=\frac{1}{x}, \quad y_1=\frac{y}{x^2}, \quad z_1=\frac{z}{x^3}.
\end{equation}
In this coordinate system, we see that this system has three accessible singular points:
\begin{equation}
(x_1,y_1,z_1)=\left\{(0,-1,2),\left(0,\frac{3}{2}(1-\sqrt{5}),9(3-\sqrt{5}) \right),\left(0,\frac{3}{2}(1+\sqrt{5}),9(3+\sqrt{5}) \right) \right\}.
\end{equation}

Around the point $(x_1,y_1,z_1)=\left(0,\frac{3}{2}(1-\sqrt{5}),9(3-\sqrt{5}) \right)$, we can rewrite the system as follows.

{\bf Step 2:} We make a change of variables.
\begin{equation}
x_2=x_1, \quad y_2=y_1-\frac{3}{2}(1-\sqrt{5}), \quad z_2=z_1-9(3-\sqrt{5}).
\end{equation}
In this coordinate system, we can rewrite the system satisfying the condition \eqref{b}:
\begin{align*}
\frac{d}{dt}\begin{pmatrix}
             x_2 \\
             y_2 \\
             z_2 
             \end{pmatrix}&=\frac{1}{x_2}\left\{\begin{pmatrix}
             \frac{3}{2}(\sqrt{5}-1) & 0 & 0  \\
             0 & 6(\sqrt{5}-1) & 1  \\
             0 & -9(\sqrt{5}-3) & \frac{9}{2}(\sqrt{5}-1)
             \end{pmatrix}\begin{pmatrix}
             x_2 \\
             y_2 \\
             z_2 
             \end{pmatrix}+\cdots\right\}.
             \end{align*}
and we can obtain the local index $\left(\frac{3}{2}(\sqrt{5}-1),3(\sqrt{5}-1),\frac{15}{2}(\sqrt{5}-1) \right)$ at the point $\{(x_2,y_2,z_2)=(0,0,0)\}$. The continued ratio of the local index at the point $\{(x_2,y_2,z_2)=(0,0,0)\}$ are all positive integers
\begin{equation}
\left(\frac{3(\sqrt{5}-1)}{\frac{3}{2}(\sqrt{5}-1)},\frac{\frac{15}{2}(\sqrt{5}-1)}{\frac{3}{2}(\sqrt{5}-1)} \right)=(2,5).
\end{equation}
This is the reason why we choose this accessible singular point.

We aim to obtain the local index $(1,0,2)$ by successive blowing-up procedures.

{\bf Step 3:} We blow up at the point $\{(x_2,y_2,z_2)=(0,0,0)\}$.
\begin{equation}
x_3=x_2, \quad y_3=\frac{y_2}{x_2}, \quad z_3=\frac{z_2}{x_2}.
\end{equation}

{\bf Step 4:} We blow up at the point $\{(x_3,y_3,z_3)=(0,0,0)\}$.
\begin{equation}
x_4=x_3, \quad y_4=\frac{y_3}{x_3}, \quad z_4=\frac{z_3}{x_3}.
\end{equation}
In this coordinate system, we see that this system has the following accessible singular locus:
\begin{equation}
(x_4,y_4,z_4)=(0,y_4,-3(\sqrt{5}-1)y_4).
\end{equation}

{\bf Step 5:} We blow up along the curve $\{(x_4,y_4,z_4)=(0,y_4,-3(\sqrt{5}-1)y_4)\}$.
\begin{equation}
x_5=x_4, \quad y_5=y_4, \quad z_5=\frac{z_4+3(\sqrt{5}-1)y_4}{x_4}.
\end{equation}

{\bf Step 6:} We make a change of variables..
\begin{equation}
x_6=\frac{1}{x_5}, \quad y_6=y_5, \quad z_6=z_5,
\end{equation}
and we can obtain the system \eqref{system;ChazyIX}.

Thus, we have completed the proof of Theorem \ref{th:7.1}. \qed

We note on the remaining accessible singular points.

Around the point $(x_1,y_1,z_1)=\left(0,\frac{3}{2}(1+\sqrt{5}),9(3+\sqrt{5}) \right)$, we can rewrite the system as follows.

{\bf Step 2:} We make a change of variables.
\begin{equation}
x_2=x_1, \quad y_2=y_1-\frac{3}{2}(1+\sqrt{5}), \quad z_2=z_1-9(3+\sqrt{5}).
\end{equation}
In this coordinate system, we can rewrite the system satisfying the condition \eqref{b}:
\begin{align*}
\frac{d}{dt}\begin{pmatrix}
             x_2 \\
             y_2 \\
             z_2 
             \end{pmatrix}&=\frac{1}{x_2}\left\{\begin{pmatrix}
             -\frac{3}{2}(\sqrt{5}+1) & 0 & 0  \\
             0 & -6(\sqrt{5}+1) & 1  \\
             0 & 9(\sqrt{5}+3) & -\frac{9}{2}(\sqrt{5}+1)
             \end{pmatrix}\begin{pmatrix}
             x_2 \\
             y_2 \\
             z_2 
             \end{pmatrix}+\cdots\right\}.
             \end{align*}
and we can obtain the local index $\left(-\frac{3}{2}(\sqrt{5}+1),-3(\sqrt{5}+1),-\frac{15}{2}(\sqrt{5}+1) \right)$ at the point $\{(x_2,y_2,z_2)=(0,0,0)\}$. The continued ratio of the local index at the point $\{(x_2,y_2,z_2)=(0,0,0)\}$ are all positive integers
\begin{equation}
\left(\frac{-3(\sqrt{5}+1)}{-\frac{3}{2}(\sqrt{5}+1)},\frac{-\frac{15}{2}(\sqrt{5}+1)}{-\frac{3}{2}(\sqrt{5}+1)} \right)=(2,5).
\end{equation}
We remark that we can discuss this case in the same way as in the case of $(x_1,y_1,z_1)=\left(0,\frac{3}{2}(1-\sqrt{5}),9(3-\sqrt{5}) \right)$.

Around the point $(x_1,y_1,z_1)=(0,-1,2)$, we can rewrite the system as follows.

{\bf Step 2:} We make a change of variables.
\begin{equation}
x_2=x_1, \quad y_2=y_1+1, \quad z_2=z_1-2.
\end{equation}
In this coordinate system, we can rewrite the system satisfying the condition \eqref{b}:
\begin{align*}
\frac{d}{dt}\begin{pmatrix}
             x_2 \\
             y_2 \\
             z_2 
             \end{pmatrix}&=\frac{1}{x_2}\left\{\begin{pmatrix}
             1 & 0 & 0  \\
             0 & 4 & 1  \\
             0 & 42 & 3
             \end{pmatrix}\begin{pmatrix}
             x_2 \\
             y_2 \\
             z_2 
             \end{pmatrix}+\cdots\right\}.
             \end{align*}
and we can obtain the local index $(1,-3,10)$ at the point $\{(x_2,y_2,z_2)=(0,0,0)\}$. The continued ratio of the local index at the point $\{(x_2,y_2,z_2)=(0,0,0)\}$ are given by
\begin{equation}
\left(\frac{-3}{1},\frac{10}{1} \right)=(-3,10).
\end{equation}
In this case, the local index involves a negative integer. So, we need to blow down.

\section{Birational symmetry of Chazy-IX equation}
In this section, let us consider the birational B{\"a}cklund transformations of the Chazy IX equation. Two birational  B{\"a}cklund transformations are new.

\begin{figure}[h]
\unitlength 0.1in
\begin{picture}( 31.7500, 12.4000)( 16.9000,-19.3000)
\put(16.9000,-9.9000){\makebox(0,0)[lb]{$(u,u',u'')$}}%
%
\special{pn 20}%
\special{pa 2066 1090}%
\special{pa 2066 1760}%
\special{fp}%
\special{sh 1}%
\special{pa 2066 1760}%
\special{pa 2086 1694}%
\special{pa 2066 1708}%
\special{pa 2046 1694}%
\special{pa 2066 1760}%
\special{fp}%
\put(16.9000,-19.7000){\makebox(0,0)[lb]{$(x,y,z)$}}%
%
\special{pn 20}%
\special{pa 2630 1870}%
\special{pa 4260 1870}%
\special{fp}%
\special{sh 1}%
\special{pa 4260 1870}%
\special{pa 4194 1850}%
\special{pa 4208 1870}%
\special{pa 4194 1890}%
\special{pa 4260 1870}%
\special{fp}%
\put(44.9000,-10.1000){\makebox(0,0)[lb]{$(p,q,r)=(p,p',p'')$}}%
\put(44.9000,-19.9000){\makebox(0,0)[lb]{$(X,Y,Z)$}}%
%
\special{pn 20}%
\special{pa 4826 1760}%
\special{pa 4826 1100}%
\special{fp}%
\special{sh 1}%
\special{pa 4826 1100}%
\special{pa 4806 1168}%
\special{pa 4826 1154}%
\special{pa 4846 1168}%
\special{pa 4826 1100}%
\special{fp}%
%
\special{pn 8}%
\special{pa 2630 900}%
\special{pa 4240 900}%
\special{fp}%
\special{sh 1}%
\special{pa 4240 900}%
\special{pa 4174 880}%
\special{pa 4188 900}%
\special{pa 4174 920}%
\special{pa 4240 900}%
\special{fp}%
\put(31.7000,-8.6000){\makebox(0,0)[lb]{$g_0$}}%
\put(20.9500,-15.1000){\makebox(0,0)[lb]{$\varphi_0$}}%
\put(48.6500,-15.1000){\makebox(0,0)[lb]{$\varphi_1$}}%
\put(31.9000,-21.0000){\makebox(0,0)[lb]{$s_0$}}%
\end{picture}%
\label{ChazyIXfig1}
\caption{}
\end{figure}

\begin{proposition}
The birational transformation $s_0$
\begin{equation}\label{bira:1}
  \left\{
  \begin{aligned}
   X &=x+\frac{2\{3(\sqrt{5}+3)y^2+\delta\}}{3(\sqrt{5}-1)z},\\
   Y &=-y,\\
   Z &=-z
   \end{aligned}
  \right. 
\end{equation}
takes the system \eqref{system;ChazyIX} to the system
\begin{equation}\label{system;1}
  \left\{
  \begin{aligned}
   \frac{dX}{dt} &=-\frac{3}{2}(\sqrt{5}-1)X^2-(9+4\sqrt{5})Y,\\
   \frac{dY}{dt} &=Z,\\
   \frac{dZ}{dt} &=3(\sqrt{5}+3)Y^2+3(\sqrt{5}-1)XZ+\delta.
   \end{aligned}
  \right. 
\end{equation}
\end{proposition}

\begin{proposition}
The birational transformation $\varphi_1$
\begin{equation}\label{bira:2}
  \left\{
  \begin{aligned}
   p &=\frac{1}{2}(\sqrt{5}-3)X,\\
   q &=3(\sqrt{5}-2)X^2+\frac{1}{2}(7+3\sqrt{5})Y,\\
   r &=9(3\sqrt{5}-7)X^3-6(\sqrt{5}+2)XY+\frac{1}{2}(3\sqrt{5}+7)Z
   \end{aligned}
  \right. 
\end{equation}
takes the system \eqref{system;1} to the system
\begin{equation}
  \left\{
  \begin{aligned}
   \frac{dp}{dt} &=q,\\
   \frac{dq}{dt} &=r,\\
   \frac{dr}{dt} &=54p^4+72p^2 q+12q^2+\frac{7+3\sqrt{5}}{2} \delta.
   \end{aligned}
  \right. 
\end{equation}
\end{proposition}
Then, we can obtain the Chazy IX equation:
\begin{equation}
\frac{d^3 p}{dt^3}=54p^4+72p^2\frac{dp}{dt}+12\left(\frac{dp}{dt} \right)^2+\frac{7+3\sqrt{5}}{2} \delta \quad (\delta \in {\Bbb C}).
\end{equation}
The compositions of the transformations \eqref{bira:1} and \eqref{bira:2} are a B{\"a}cklund transformation of the Chazy IX equation (see Figure 2).
\begin{theorem}
The Chazy IX equation \eqref{eq;ChazyIX} is invariant under the biratioral transformation$:$
\begin{align}
\begin{split}
&g_0(u;\delta) \rightarrow \\
&\left(\frac{(\sqrt{5}-3)\{108u^4+18(5+\sqrt{5})u^2 u'+6(3+\sqrt{5})(u')^2+3(\sqrt{5}-1)uu''+2\delta}{6(\sqrt{5}-1)\{3(\sqrt{5}-1)uu'+u''\}};\frac{7+3\sqrt{5}}{2} \delta \right),
\end{split}
\end{align}
where $u'=\frac{du}{dt},\ u''=\frac{d^2u}{dt^2}$.
\end{theorem}
This B{\"a}cklund transformation is new.

\begin{figure}[ht]
\unitlength 0.1in
\begin{picture}( 58.5500, 13.0000)( 16.9000,-19.5000)
\put(16.9000,-9.9000){\makebox(0,0)[lb]{$(u,u',u'')$}}%
%
\special{pn 20}%
\special{pa 2066 1090}%
\special{pa 2066 1760}%
\special{fp}%
\special{sh 1}%
\special{pa 2066 1760}%
\special{pa 2086 1694}%
\special{pa 2066 1708}%
\special{pa 2046 1694}%
\special{pa 2066 1760}%
\special{fp}%
\put(16.9000,-19.7000){\makebox(0,0)[lb]{$(x,y,z)$}}%
%
\special{pn 20}%
\special{pa 2630 1870}%
\special{pa 4260 1870}%
\special{fp}%
\special{sh 1}%
\special{pa 4260 1870}%
\special{pa 4194 1850}%
\special{pa 4208 1870}%
\special{pa 4194 1890}%
\special{pa 4260 1870}%
\special{fp}%
\put(44.9000,-19.9000){\makebox(0,0)[lb]{$(X,Y,Z)$}}%
\put(43.0000,-8.6000){\makebox(0,0)[lb]{$g_1$}}%
\put(20.9500,-15.1000){\makebox(0,0)[lb]{$\varphi_0$}}%
\put(31.9000,-21.0000){\makebox(0,0)[lb]{$\pi$}}%
%
\special{pn 20}%
\special{pa 5310 1890}%
\special{pa 6940 1890}%
\special{fp}%
\special{sh 1}%
\special{pa 6940 1890}%
\special{pa 6874 1870}%
\special{pa 6888 1890}%
\special{pa 6874 1910}%
\special{pa 6940 1890}%
\special{fp}%
\put(71.7000,-8.2000){\makebox(0,0)[lb]{$(p,q,r)$}}%
\put(71.7000,-20.1000){\makebox(0,0)[lb]{$(u,v,w)$}}%
%
\special{pn 20}%
\special{pa 7506 1780}%
\special{pa 7506 1120}%
\special{fp}%
\special{sh 1}%
\special{pa 7506 1120}%
\special{pa 7486 1188}%
\special{pa 7506 1174}%
\special{pa 7526 1188}%
\special{pa 7506 1120}%
\special{fp}%
\put(75.4500,-15.3000){\makebox(0,0)[lb]{$\varphi_2$}}%
\put(58.7000,-21.2000){\makebox(0,0)[lb]{$s_1$}}%
%
\special{pn 8}%
\special{pa 2630 910}%
\special{pa 6920 910}%
\special{fp}%
\special{sh 1}%
\special{pa 6920 910}%
\special{pa 6854 890}%
\special{pa 6868 910}%
\special{pa 6854 930}%
\special{pa 6920 910}%
\special{fp}%
\put(70.4000,-10.1000){\makebox(0,0)[lb]{$=(p,p',p'')$}}%
\end{picture}%
\label{ChazyIXfig2}
\caption{}
\end{figure}

\begin{proposition}
The birational transformation $\pi$
\begin{equation}\label{bira:0}
  \left\{
  \begin{aligned}
   X &=x,\\
   Y &=y-3\sqrt{5}x^2,\\
   Z &=z-6\sqrt{5}xy-9(\sqrt{5}-5)x^3
   \end{aligned}
  \right. 
\end{equation}
takes the system \eqref{system;ChazyIX} to the system
\begin{equation}\label{system;ChazyIXb}
  \left\{
  \begin{aligned}
   \frac{dX}{dt} &=-\frac{3}{2}(-\sqrt{5}-1)X^2+Y,\\
   \frac{dY}{dt} &=Z,\\
   \frac{dZ}{dt} &=3(-\sqrt{5}+3)Y^2+3(-\sqrt{5}-1)XZ+\delta.
   \end{aligned}
  \right. 
\end{equation}
\end{proposition}
This transformation changes the sign of $\sqrt{5}$ in the system \eqref{system;ChazyIX} (cf. \cite{Cos1}).

\begin{proposition}
The birational transformation $s_1$
\begin{equation}\label{bira:3}
  \left\{
  \begin{aligned}
   u &=X-\frac{2\{3(-\sqrt{5}+3)Y^2+\delta\}}{3(\sqrt{5}+1)Z},\\
   v &=-Y,\\
   w &=-Z
   \end{aligned}
  \right. 
\end{equation}
takes the system \eqref{system;ChazyIXb} to the system
\begin{equation}\label{system;2}
  \left\{
  \begin{aligned}
   \frac{du}{dt} &=-\frac{3}{2}(-\sqrt{5}-1)u^2+(-9+4\sqrt{5})v,\\
   \frac{dv}{dt} &=w,\\
   \frac{dw}{dt} &=3(-\sqrt{5}+3)v^2+3(-\sqrt{5}-1)uw+\delta.
   \end{aligned}
  \right. 
\end{equation}
\end{proposition}

\begin{proposition}
The birational transformation $\varphi_2$
\begin{equation}\label{bira:4}
  \left\{
  \begin{aligned}
   p &=\frac{1}{2}(-\sqrt{5}-3)u,\\
   q &=3(-\sqrt{5}-2)u^2+\frac{1}{2}(7-3\sqrt{5})v,\\
   r &=9(-3\sqrt{5}-7)u^3-6(-\sqrt{5}+2)uv+\frac{1}{2}(-3\sqrt{5}+7)w
   \end{aligned}
  \right. 
\end{equation}
takes the system \eqref{system;2} to the system
\begin{equation}
  \left\{
  \begin{aligned}
   \frac{dp}{dt} &=q,\\
   \frac{dq}{dt} &=r,\\
   \frac{dr}{dt} &=54p^4+72p^2 q+12q^2+\frac{7-3\sqrt{5}}{2} \delta.
   \end{aligned}
  \right. 
\end{equation}
\end{proposition}
Then, we can obtain the Chazy IX equation:
\begin{equation}
\frac{d^3 p}{dt^3}=54p^4+72p^2\frac{dp}{dt}+12\left(\frac{dp}{dt} \right)^2+\frac{7-3\sqrt{5}}{2} \delta \quad (\delta \in {\Bbb C}).
\end{equation}
The compositions of the transformations \eqref{bira:0}, \eqref{bira:3} and \eqref{bira:4} are a B{\"a}cklund transformation of the Chazy IX equation (see Figure 3).
\begin{theorem}
The Chazy IX equation \eqref{eq;ChazyIX} is invariant under the biratioral transformation$:$
\begin{align}
\begin{split}
&g_1(u;\delta) \rightarrow \\
&\left(\frac{(-\sqrt{5}-3)\{108u^4+18(5-\sqrt{5})u^2 u'+6(3-\sqrt{5})(u')^2+3(-\sqrt{5}-1)uu''+2\delta}{6(-\sqrt{5}-1)\{3(-\sqrt{5}-1)uu'+u''\}};\frac{7-3\sqrt{5}}{2} \delta \right),
\end{split}
\end{align}
where $u'=\frac{du}{dt},\ u''=\frac{d^2u}{dt^2}$.
\end{theorem}
This B{\"a}cklund transformation is new.

\section{Holomorphy conditions of Chazy-IX equation}
In this section, we give the holomorphy condition of the system \eqref{system;ChazyIX}. Thanks to this holomorphy condition, we obtain a new partial differential system in two variables involving the Chazy IX equation, This system satisfies the compatibility condition, and admits a travelling wave solution.
\begin{theorem}
Let us consider the following differential system in the polynomial class\rm{:\rm}
\begin{equation*}
  \left\{
  \begin{aligned}
   dx &=f_1(x,y,z)dt+g_1(x,y,z)ds,\\
   dy &=f_2(x,y,z)dt+g_2(x,y,z)ds,\\
   dz &=f_3(x,y,z)dt+g_3(x,y,z)ds.
   \end{aligned}
  \right. 
\end{equation*}
We assume that

$(A1)$ $deg(f_i)=3$ and $deg(g_i)=3$ with respect to $x,y,z$.

$(A2)$ The right-hand side of this system becomes again a polynomial in each coordinate system $(x_i,y_i,z_i) \ (i=1,2)$.
\begin{align}
\begin{split}
1) \ &x_1=\frac{1}{x}, \quad y_1=y, \quad z_1=-\left(zx+\frac{2\{3(\sqrt{5}+3)y^2+\delta\}}{3(\sqrt{5}-1)} \right)x,\\
2) \ &x_2=\frac{1}{x}, \quad y_2=y-3\sqrt{5} x^2, \quad z_2=-((z-6\sqrt{5}xy-9(\sqrt{5}-5)x^3)x\\
&-\frac{2\{-(\sqrt{5}-3)(135x^4+3y^2)+(90-54\sqrt{5})x^2y+\delta\}}{3(\sqrt{5}+1)})x.
\end{split}
\end{align}
Then such a system coincides with
\begin{equation}\label{total}
  \left\{
  \begin{aligned}
   dx =&\left\{-\frac{3}{2}(\sqrt{5}-1)x^2+y\right\}dt\\
   &+\left\{-12(\sqrt{5}+2)x^2y+(\sqrt{5}+1)\left(3xz-2y^2-\frac{2}{3}\delta \right) \right\}ds,\\
   dy =&z dt\\
   &+\{-6(2\sqrt{5}-5)x^2z+12\sqrt{5}xy^2-(\sqrt{5}+1)yz+(3\sqrt{5}-5)\delta x \}ds,\\
   dz =&\{3(\sqrt{5}+3)y^2+3(\sqrt{5}-1)xz+\delta\}dt\\
   &+\{48xyz-24y^3-(\sqrt{5}+1)z^2+2(\sqrt{5}-3)\delta y \}ds.
   \end{aligned}
  \right. 
\end{equation}
\end{theorem}
These transition functions satisfy the condition{\rm:\rm}
\begin{equation*}
dx_i \wedge dy_i \wedge dz_i=dx \wedge dy \wedge dz \quad (i=1,2).
\end{equation*}

We remark that when $s=0$, we can obtain the system \eqref{system;ChazyIX}.

\begin{proposition}
The system \eqref{total} satisfies the compatibility condition:
\begin{equation}
\frac{\partial }{\partial s} \frac{\partial x}{\partial t}=\frac{\partial }{\partial t} \frac{\partial x}{\partial s}, \quad \frac{\partial }{\partial s} \frac{\partial y}{\partial t}=\frac{\partial }{\partial t} \frac{\partial y}{\partial s}, \quad \frac{\partial }{\partial s} \frac{\partial z}{\partial t}=\frac{\partial }{\partial t} \frac{\partial z}{\partial s}.
\end{equation}
\end{proposition}

\begin{proposition}
The system \eqref{total} admits the travelling wave solution{\rm:\rm}
\begin{align}
\begin{split}
(x,y,z;\delta)=&(\frac{c_2}{\sqrt{3(\sqrt{5}+1)}} Tanh \left(\sqrt{\frac{3}{2}(\sqrt{5}-2)} c_2 t+\sqrt{3(\sqrt{5}+1)} c_2^3 s+c_1 \right),\\
&-\frac{1}{4}(\sqrt{5}-3)c_2^2,0;-\frac{3}{8}(\sqrt{5}-3)c_2^4 ) \quad (c_i \in {\Bbb C}).
\end{split}
\end{align}
\end{proposition}

\section{Chazy IX equation and soliton equations}
In this section, we consider the relation between the system \eqref{total} and soliton equations. In this paper, we can make the birational transformations between the system \eqref{total} and soliton equations.
\begin{theorem}
The birational transformations
\begin{equation}\label{eq:A2}
  \left\{
  \begin{aligned}
   X =&x,\\
   Y =&y+\frac{3(1-\sqrt{5})}{2}x^2,\\
   Z =&z+3(1-\sqrt{5})xy+9(3-\sqrt{5})x^3,\\
   S =&\frac{8}{3(1-\sqrt{5})}s
   \end{aligned}
  \right. 
\end{equation}
take the system \eqref{total} to the system
\begin{equation}\label{eq:A3}
  \left\{
  \begin{aligned}
   dX =&Y dt+\left(54X^4+18X^2Y+3Y^2-\frac{9}{2}XZ+\delta \right)dS,\\
   dY =&Z dt+\left(-243X^5-108X^3 Y-18X Y^2+18X^2 Z+\frac{3}{2}YZ-\frac{9}{2}\delta X \right)dS,\\
   dZ =&(54X^4+72X^2Y+12Y^2+\delta) dt\\
&+\left(972X^6+162X^4Y-108X^3Z+\frac{3}{2}Z^2+18 \delta X^2-3 \delta Y \right)dS.
   \end{aligned}
  \right. 
\end{equation}
\end{theorem}
Setting $u:=X$, we see that
\begin{equation}
\frac{\partial u}{\partial t}=Y, \quad  \frac{\partial^2 u}{\partial t^2}=Z,
\end{equation}
and
\begin{equation}\label{eq:A4}
  \left\{
  \begin{aligned}
   \frac{\partial^3 u}{\partial t^3} =&54u^4+72u^2 \frac{\partial u}{\partial t}+12\left( \frac{\partial u}{\partial t} \right)^2+\delta,\\
   \frac{\partial u}{\partial S} =&54u^4+18u^2 \frac{\partial u}{\partial t}+3\left( \frac{\partial u}{\partial t} \right)^2-\frac{9}{2}u\frac{\partial^2 u}{\partial t^2}+\delta.
   \end{aligned}
  \right. 
\end{equation}
The first equation in \eqref{eq:A4} coincides with Chazy IX equation.

Adding each system in \eqref{eq:A4}, we can obtain
\begin{equation}\label{eq:A5}
   \frac{\partial^3 u}{\partial t^3} =54u^2\frac{\partial u}{\partial t}+9\left( \frac{\partial u}{\partial t} \right)^2+\frac{9}{2}u \frac{\partial^2 u}{\partial t^2}+\frac{\partial u}{\partial S}.
\end{equation}

\begin{question}
It is still an open question whether the second equation in \eqref{eq:A4} and the system \eqref{eq:A5} coincide with which of the soliton equations.
\end{question}

\section{The Chazy X equation}
The Chazy-X.a equation is given by
\begin{align}\label{eq;ChazyX}
\begin{split}
\frac{d^3 u}{dt^3}=&6u^2 \frac{du}{dt}+\frac{3}{11}(9+7\sqrt{3})\left(\frac{du}{dt}+u^2 \right)^2\\
&-\frac{1}{22}(4-3\sqrt{3})\alpha \frac{du}{dt}+\frac{1}{44}(3-5\sqrt{3})\alpha u^2-\frac{1}{352}(9+7\sqrt{3}){\alpha}^2.
\end{split}
\end{align}
In this section, at first we transform the equation \eqref{eq;ChazyX} to a system of differential equations by birational transformations. For this system, we give two B{\"a}cklund transformations from Chazy-X.a to Chazy-X.b. One of them is new. We also give the holomorphy condition of this system. Thanks to this condition, we can recover the Chazy-X.a system.

\begin{theorem}\label{th:11.1}
The birational transformation $\varphi_0$
\begin{equation}
  \left\{
  \begin{aligned}
   x &=u,\\
   y &=\frac{du}{dt}-\frac{3+\sqrt{3}}{2}u^2,\\
   z &=\frac{d^2u}{dt^2}-(3+\sqrt{3})u \frac{du}{dt}
   \end{aligned}
  \right. 
\end{equation}
takes the Chazy X.a equation to the system of the first-order ordinary differential equations:
\begin{equation}\label{system;ChazyX}
  \left\{
  \begin{aligned}
   \frac{dx}{dt} &=\frac{3+\sqrt{3}}{2}x^2+y,\\
   \frac{dy}{dt} &=z,\\
   \frac{dz}{dt} &=\frac{2}{11}(-3+5\sqrt{3})y^2-(3+\sqrt{3})xz+\frac{1}{22}(-4+3\sqrt{3})\alpha y-\frac{1}{352}(9+7\sqrt{3})\alpha^2.
   \end{aligned}
  \right. 
\end{equation}
\end{theorem}

By the same way of the second Painlev\'e system, we can prove Theorem \ref{th:11.1}.

{\bf Proof.} At first, we rewrite the equation \eqref{eq;ChazyX} to the system of the first-order ordinary differential equations.

{\bf Step 0:} We make a change of variables.
\begin{equation}
x=u, \quad y=\frac{du}{dt}, \quad z=\frac{d^2u}{dt^2}.
\end{equation}

{\bf Step 1:} We make a change of variables.
\begin{equation}
x_1=\frac{1}{x}, \quad y_1=\frac{y}{x^2}, \quad z_1=\frac{z}{x^3}.
\end{equation}
In this coordinate system, we see that this system has three accessible singular points:
\begin{equation}
(x_1,y_1,z_1)=\left\{(0,-1,2),\left(0,\frac{3+\sqrt{3}}{2},3(2+\sqrt{3}) \right),\left(0,\frac{-1-2\sqrt{3}}{11},\frac{2}{121}(13+4\sqrt{3}) \right) \right\}.
\end{equation}

Around the point $(x_1,y_1,z_1)=\left(0,\frac{3+\sqrt{3}}{2},3(2+\sqrt{3}) \right)$, we can rewrite the system as follows.

{\bf Step 2:} We make a change of variables.
\begin{equation}
x_2=x_1, \quad y_2=y_1-\frac{3+\sqrt{3}}{2}, \quad z_2=z_1-3(2+\sqrt{3}).
\end{equation}
In this coordinate system, we can rewrite the system satisfying the condition \eqref{b}:
\begin{align*}
\frac{d}{dt}\begin{pmatrix}
             x_2 \\
             y_2 \\
             z_2 
             \end{pmatrix}&=\frac{1}{x_2}\left\{\begin{pmatrix}
             -\frac{3+\sqrt{3}}{2} & 0 & 0  \\
             0 & -2(3+\sqrt{3}) & 1  \\
             0 & 3(2+\sqrt{3}) & -\frac{3(3+\sqrt{3})}{2}
             \end{pmatrix}\begin{pmatrix}
             x_2 \\
             y_2 \\
             z_2 
             \end{pmatrix}+\cdots\right\}.
             \end{align*}
and we can obtain the local index $\left(-\frac{3+\sqrt{3}}{2},-(3+\sqrt{3}),-\frac{5}{2}(3+\sqrt{3}) \right)$ at the point $\{(x_2,y_2,z_2)=(0,0,0)\}$. The continued ratio of the local index at the point $\{(x_2,y_2,z_2)=(0,0,0)\}$ are all positive integers
\begin{equation}
\left(\frac{-(3+\sqrt{3})}{-\frac{3+\sqrt{3}}{2}},\frac{-\frac{5}{2}(3+\sqrt{3})}{-\frac{3+\sqrt{3}}{2}} \right)=(2,5).
\end{equation}
This is the reason why we choose this accessible singular point.

We aim to obtain the local index $(1,0,2)$ by successive blowing-up procedures.

{\bf Step 3:} We blow up at the point $\{(x_2,y_2,z_2)=(0,0,0)\}$.
\begin{equation}
x_3=x_2, \quad y_3=\frac{y_2}{x_2}, \quad z_3=\frac{z_2}{x_2}.
\end{equation}

{\bf Step 4:} We blow up at the point $\{(x_3,y_3,z_3)=(0,0,0)\}$.
\begin{equation}
x_4=x_3, \quad y_4=\frac{y_3}{x_3}, \quad z_4=\frac{z_3}{x_3}.
\end{equation}
In this coordinate system, we see that this system has the following accessible singular locus:
\begin{equation}
(x_4,y_4,z_4)=(0,y_4,(3+\sqrt{3})y_4).
\end{equation}

{\bf Step 5:} We blow up along the curve $\{(x_4,y_4,z_4)=(0,y_4,(3+\sqrt{3})y_4)\}$.
\begin{equation}
x_5=x_4, \quad y_5=y_4, \quad z_5=\frac{z_4-(3+\sqrt{3})y_4}{x_4}.
\end{equation}

{\bf Step 6:} We make a change of variables..
\begin{equation}
x_6=\frac{1}{x_5}, \quad y_6=y_5, \quad z_6=z_5,
\end{equation}
and we can obtain the system \eqref{system;ChazyX}.

Thus, we have completed the proof of Theorem \ref{th:11.1}. \qed

We note on the remaining accessible singular points.

Around the point $(x_1,y_1,z_1)=(0,-1,2)$, we can rewrite the system as follows.

{\bf Step 2:} We make a change of variables.
\begin{equation}
x_2=x_1, \quad y_2=y_1+1, \quad z_2=z_1-2.
\end{equation}
In this coordinate system, we can rewrite the system satisfying the condition \eqref{b}:
\begin{align*}
\frac{d}{dt}\begin{pmatrix}
             x_2 \\
             y_2 \\
             z_2 
             \end{pmatrix}&=\frac{1}{x_2}\left\{\begin{pmatrix}
             1 & 0 & 0  \\
             0 & 4 & 1  \\
             0 & 0 & 3
             \end{pmatrix}\begin{pmatrix}
             x_2 \\
             y_2 \\
             z_2 
             \end{pmatrix}+\cdots\right\}.
             \end{align*}
and we can obtain the local index $(1,4,3)$ at the point $\{(x_2,y_2,z_2)=(0,0,0)\}$. The continued ratio of the local index at the point $\{(x_2,y_2,z_2)=(0,0,0)\}$ are all positive integers
\begin{equation}
\left(\frac{4}{1},\frac{3}{1} \right)=(4,3).
\end{equation}
We remark that we can discuss this case in the same way as in the case of $(x_1,y_1,z_1)=\left(0,\frac{3+\sqrt{3}}{2},3(2+\sqrt{3}) \right)$.

Around the point $(x_1,y_1,z_1)=\left(0,\frac{-1-2\sqrt{3}}{11},\frac{2}{121}(13+4\sqrt{3}) \right)$, we can rewrite the system as follows.

{\bf Step 2:} We make a change of variables.
\begin{equation}
x_2=x_1, \quad y_2=y_1+\frac{1+2\sqrt{3}}{11}, \quad z_2=z_1-\frac{2}{121}(13+4\sqrt{3}).
\end{equation}
In this coordinate system, we can rewrite the system satisfying the condition \eqref{b}:
\begin{align*}
\frac{d}{dt}\begin{pmatrix}
             x_2 \\
             y_2 \\
             z_2 
             \end{pmatrix}&=\frac{1}{x_2}\left\{\begin{pmatrix}
             \frac{1+2\sqrt{3}}{11} & 0 & 0  \\
             0 & \frac{4(1+2\sqrt{3})}{11} & 1  \\
             0 & \frac{72(13+4\sqrt{3})}{121} & \frac{3(1+2\sqrt{3})}{11}
             \end{pmatrix}\begin{pmatrix}
             x_2 \\
             y_2 \\
             z_2 
             \end{pmatrix}+\cdots\right\}.
             \end{align*}
and we can obtain the local index $\left(\frac{1+2\sqrt{3}}{11},\frac{12(1+2\sqrt{3})}{11},\frac{-5(1+2\sqrt{3})}{11} \right)$ at the point $\{(x_2,y_2,z_2)=(0,0,0)\}$. The continued ratio of the local index at the point $\{(x_2,y_2,z_2)=(0,0,0)\}$ are given by
\begin{equation}
\left(\frac{\frac{12(1+2\sqrt{3})}{11}}{\frac{1+2\sqrt{3}}{11}},\frac{\frac{-5(1+2\sqrt{3})}{11}}{\frac{1+2\sqrt{3}}{11}} \right)=(12,-5).
\end{equation}
In this case, the local index involves a negative integer. So, we need to blow down.

\section{Birational symmetry of Chazy X equation}
In this section, let us consider the birational B{\"a}cklund transformations from Chazy-X.a to Chazy-X.b. One of them is new.

\begin{proposition}
The birational transformation $s_0$
\begin{equation}\label{bira:11}
  \left\{
  \begin{aligned}
   X &=x-\frac{64(-3+5\sqrt{3})y^2+16(-4+3\sqrt{3})\alpha y-(9+7\sqrt{3})\alpha^2}{176(3+\sqrt{3})z},\\
   Y &=-y,\\
   Z &=-z
   \end{aligned}
  \right. 
\end{equation}
takes the system \eqref{system;ChazyX} to the system
\begin{equation}\label{system;11}
  \left\{
  \begin{aligned}
   \frac{dX}{dt} &=\frac{3+\sqrt{3}}{2}X^2+\frac{1}{11}(-43+24\sqrt{3})Y+\frac{21-13\sqrt{3}}{66}\alpha,\\
   \frac{dY}{dt} &=Z,\\
   \frac{dZ}{dt} &=\frac{2}{11}(-3+5\sqrt{3})Y^2-(3+\sqrt{3})XZ-\frac{1}{22}(-4+3\sqrt{3})\alpha Y-\frac{1}{352}(9+7\sqrt{3})\alpha^2.
   \end{aligned}
  \right. 
\end{equation}
\end{proposition}

\begin{proposition}
The birational transformation $\varphi_1$
\begin{equation}\label{bira:12}
  \left\{
  \begin{aligned}
   p &=(2+\sqrt{3})X,\\
   q &=\frac{2+\sqrt{3}}{11(3+\sqrt{3})}\{33(2+\sqrt{3})X^2+(-57+29\sqrt{3})Y-(-4+3\sqrt{3})\alpha\},\\
   r &=\frac{6(2+\sqrt{3})}{11(3+\sqrt{3})^3} \{33(33+19\sqrt{3})X^3-6(13+4\sqrt{3})XY+(-27+\sqrt{3})Z-(9+7\sqrt{3})\alpha X\}
   \end{aligned}
  \right. 
\end{equation}
takes the system \eqref{system;11} to the system
\begin{equation}
  \left\{
  \begin{aligned}
   \frac{dp}{dt} &=q,\\
   \frac{dq}{dt} &=r,\\
   \frac{dr}{dt} &=6p^2q+\frac{3}{11}(9-7\sqrt{3})\left(q+p^2 \right)^2\\
&-\frac{1}{22}(4+3\sqrt{3})\alpha q+\frac{1}{44}(3+5\sqrt{3})\alpha p^2-\frac{1}{352}(9-7\sqrt{3}){\alpha}^2.
   \end{aligned}
  \right. 
\end{equation}
\end{proposition}
Then, we can obtain the Chazy X.b equation:
\begin{align}
\begin{split}
\frac{d^3 u}{dt^3}=&6u^2 \frac{du}{dt}+\frac{3}{11}(9-7\sqrt{3})\left(\frac{du}{dt}+u^2 \right)^2\\
&-\frac{1}{22}(4+3\sqrt{3})\alpha \frac{du}{dt}+\frac{1}{44}(3+5\sqrt{3})\alpha u^2-\frac{1}{352}(9-7\sqrt{3}){\alpha}^2,
\end{split}
\end{align}
The compositions of the transformations \eqref{bira:11} and \eqref{bira:12} are a B{\"a}cklund transformation from Chazy-X.a to Chazy-X.b (see Figure 1).
\begin{theorem}
The Chazy X.a equation \eqref{eq;ChazyX} can be transformed into the Chazy X.b equation by the birational transformation$:$
\begin{align}
\begin{split}
&g_0(u;\alpha) \rightarrow \\
&(\frac{3+\sqrt{3}}{1056\{(3+\sqrt{3})uu'-u''\}} \{288(6+\sqrt{3})u^2 u'-176(3+\sqrt{3})uu''+96(9+7\sqrt{3})u^4\\
&+64(-3+5\sqrt{3})(u')^2-8(-3+5\sqrt{3})\alpha u^2+16(-4+3\sqrt{3})\alpha u'-(9+7\sqrt{3})\alpha^2 \};\alpha ),
\end{split}
\end{align}
where $u'=\frac{du}{dt},\ u''=\frac{d^2u}{dt^2}$.
\end{theorem}

\begin{proposition}
The birational transformation $\pi$
\begin{equation}\label{bira:10}
  \left\{
  \begin{aligned}
   X &=x,\\
   Y &=y+\frac{1}{2}(5+\sqrt{3})x^2+\frac{1}{8}(\sqrt{3}-1)\alpha,\\
   Z &=z-(5+\sqrt{3})x^3-\frac{1}{4}(2\sqrt{3}-1)\alpha x
   \end{aligned}
  \right. 
\end{equation}
takes the system \eqref{system;ChazyX} to the system
\begin{equation}\label{system;ChazyXb}
  \left\{
  \begin{aligned}
   \frac{dX}{dt} &=-X^2+Y-\frac{1}{8}(\sqrt{3}-1)\alpha,\\
   \frac{dY}{dt} &=(5+\sqrt{3})XY+Z,\\
   \frac{dZ}{dt} &=-(15+7\sqrt{3})X^2 Y+\frac{2}{11}(-3+5\sqrt{3})Y^2-(3+\sqrt{3})XZ-\frac{3}{4}\alpha Y.
   \end{aligned}
  \right. 
\end{equation}
\end{proposition}

\begin{proposition}
The birational transformation $s_1$
\begin{equation}\label{bira:13}
  \left\{
  \begin{aligned}
   u &=X-\frac{(\sqrt{3}-4)Z}{13Y},\\
   v &=Y,\\
   w &=-Z
   \end{aligned}
  \right. 
\end{equation}
takes the system \eqref{system;ChazyXb} to the system
\begin{equation}\label{system;12}
  \left\{
  \begin{aligned}
   \frac{du}{dt} &=-(4+\sqrt{3})u^2+\frac{8(89+46\sqrt{3})v-11(11+7\sqrt{3})\alpha}{1144},\\
   \frac{dv}{dt} &=(5+\sqrt{3})uv-\frac{1}{13}(-4+\sqrt{3})w,\\
   \frac{dw}{dt} &=(15+7\sqrt{3})u^2 v+(3+\sqrt{3})uw+\frac{1}{44}v\{8(3-5\sqrt{3})v+33\alpha\}.
   \end{aligned}
  \right. 
\end{equation}
\end{proposition}

\begin{proposition}
The birational transformation $\varphi_2$
\begin{equation}\label{bira:14}
  \left\{
  \begin{aligned}
   p &=(4+\sqrt{3})u,\\
   q &=-(19+8\sqrt{3})u^2+\frac{1}{11}(38+21\sqrt{3})v-\frac{1}{8}(5+3\sqrt{3})\alpha,\\
   r &=2(100+51\sqrt{3})u^3+\frac{1}{143}(89+46\sqrt{3})w+\frac{u}{44}\{-4(177+101\sqrt{3})v+11(29+17\sqrt{3})\alpha \}
   \end{aligned}
  \right. 
\end{equation}
takes the system \eqref{system;12} to the system
\begin{equation}
  \left\{
  \begin{aligned}
   \frac{dp}{dt} &=q,\\
   \frac{dq}{dt} &=r,\\
   \frac{dr}{dt} &=6p^2q+\frac{3}{11}(9-7\sqrt{3})\left(q+p^2 \right)^2\\
&-\frac{1}{22}(4+3\sqrt{3})\beta q+\frac{1}{44}(3+5\sqrt{3})\beta p^2-\frac{1}{352}(9-7\sqrt{3}){\beta}^2,
   \end{aligned}
  \right. 
\end{equation}
where $\beta=(-2+\sqrt{3})\alpha$.
\end{proposition}
Then, we can obtain the Chazy X.b equation:
\begin{align*}
\frac{d^3 u}{dt^3}=&6u^2 \frac{du}{dt}+\frac{3}{11}(9-7\sqrt{3})\left(\frac{du}{dt}+u^2 \right)^2\\
&-\frac{1}{22}(4+3\sqrt{3})\beta \frac{du}{dt}+\frac{1}{44}(3+5\sqrt{3})\beta u^2-\frac{1}{352}(9-7\sqrt{3}){\beta}^2, \quad \beta=(-2+\sqrt{3})\alpha.
\end{align*}
The compositions of the transformations \eqref{bira:10}, \eqref{bira:13} and \eqref{bira:14} are a B{\"a}cklund transformation  from Chazy-X.a to Chazy-X.b (see Figure 2).
\begin{theorem}
The Chazy X.a equation \eqref{eq;ChazyX} can be transformed into the Chazy X.b equation by the birational transformation$:$
\begin{align}
\begin{split}
&g_1(u;\alpha) \rightarrow \\
&\left(\frac{(4+\sqrt{3})\{8(-4+\sqrt{3})u^3-8(-4+\sqrt{3})u u'-8(-4+\sqrt{3})u''-(-7+5\sqrt{3})\alpha u \}}{13\{8u^2+8u'+(\sqrt{3}-1)\alpha\}};(-2+\sqrt{3})\alpha \right),
\end{split}
\end{align}
where $u'=\frac{du}{dt},\ u''=\frac{d^2u}{dt^2}$.
\end{theorem}
This B{\"a}cklund transformation from Chazy-X.a to Chazy-X.b is new.

\section{Holomorphy conditions of Chazy X equation}
In this section, we give the holomorphy condition of the system \eqref{system;ChazyX}. Thanks to this holomorphy condition, we can recover the system \eqref{system;ChazyX}.
\begin{theorem}
Let us consider the following ordinary differential system in the polynomial class\rm{:\rm}
\begin{equation*}
  \left\{
  \begin{aligned}
   \frac{dx}{dt} &=f_1(x,y,z),\\
   \frac{dy}{dt} &=f_2(x,y,z),\\
   \frac{dz}{dt} &=f_3(x,y,z).
   \end{aligned}
  \right. 
\end{equation*}
We assume that

$(A1)$ $deg(f_i)=3$ with respect to $x,y,z$.

$(A2)$ The right-hand side of this system becomes again a polynomial in each coordinate system $(x_i,y_i,z_i) \ (i=1,2)$.
\begin{align}
\begin{split}
1) \ &x_1=\frac{1}{x}, \quad y_1=y,\\
&z_1=-\left(zx-\left(\frac{4}{11}(-4+3\sqrt{3})y^2+\frac{1}{66}(-21+13\sqrt{3})\alpha y-\frac{1+2\sqrt{3}}{176}\alpha^2 \right) \right)x,\\
2) \ &x_2=\frac{1}{x},\\
&y_2=-\left( \left(y+\frac{5+\sqrt{3}}{2}x^2+\frac{\sqrt{3}-1}{8}\alpha \right)x+\frac{-17+\sqrt{3}}{13}x^3-\frac{-4+\sqrt{3}}{13}z-\frac{-10+9\sqrt{3}}{52}\alpha x \right)x,\\
&z_2=z-(5+\sqrt{3})x^3-\frac{1}{4}(2\sqrt{3}-1)\alpha x.
\end{split}
\end{align}
Then such a system coincides with
\begin{equation}
  \left\{
  \begin{aligned}
   \frac{dx}{dt} &=-\frac{a(t)}{3+\sqrt{3}}\left(\frac{3+\sqrt{3}}{2}x^2+y \right),\\
   \frac{dy}{dt} &=-\frac{a(t)}{3+\sqrt{3}}z,\\
   \frac{dz}{dt} &=-\frac{a(t)}{3+\sqrt{3}}\left(\frac{2}{11}(-3+5\sqrt{3})y^2-(3+\sqrt{3})xz+\frac{1}{22}(-4+3\sqrt{3})\alpha y-\frac{1}{352}(9+7\sqrt{3})\alpha^2 \right),
   \end{aligned}
  \right. 
\end{equation}
where $a(t) \in {\Bbb C}(t)$. Setting $a(t)=-(3+\sqrt{3})$, we obtain the system \eqref{system;ChazyX}.
\end{theorem}
These transition functions satisfy the condition{\rm:\rm}
\begin{equation*}
dx_i \wedge dy_i \wedge dz_i=dx \wedge dy \wedge dz \quad (i=1,2).
\end{equation*}

\section{Particular solutions of Chazy X equation}
In this section, we study a solution of the system \eqref{system;ChazyX} which is written by the use of known functions.

\begin{proposition}
The system \eqref{system;ChazyX} admits two rational solutions{\rm:\rm}
\begin{align}
\begin{split}
&(x,y,z;\alpha)=(0,0,0;0)\\
&(x,y,z;\alpha)=\left(-\frac{2}{(3+\sqrt{3})t+2c},0,0;0 \right),
\end{split}
\end{align}
and  admits two special solutions{\rm:\rm}
\begin{align}
\begin{split}
(x,y,z;\alpha)=&\left(-\frac{(-3+\sqrt{3})\sqrt{3+\sqrt{3}}\sqrt{\alpha}Tan\left\{\frac{3^{\frac{1}{4}}}{4}\sqrt{3+\sqrt{3}}\sqrt{\alpha}(t+8c) \right\}}{4 \times 3^{\frac{3}{4}}},\frac{\sqrt{3}}{8}\alpha,0;\alpha \right),\\
(x,y,z;\alpha)=&\left(-\frac{(3+\sqrt{3})^{\frac{3}{2}}\sqrt{\alpha}Tanh\left\{\frac{1}{4}\sqrt{5+3\sqrt{3}}\sqrt{\alpha}(t+24c) \right\}}{12\sqrt{3+2\sqrt{3}}},-\frac{3+2\sqrt{3}}{24}\alpha,0;\alpha \right),
\end{split}
\end{align}
where $c \in {\Bbb C}$ is an integral constant.
\end{proposition}

\section{Appendix A}
In this section, we study the Chazy I equation (see \cite{Cos1}):
\begin{align}
\begin{split}
I:\frac{d^3 u}{dt^3}=&6\left\{-\left(\frac{du}{dt} \right)^2+A(t)\left(\frac{du}{dt}+u^2 \right)+B(t)u+C(t)\right\},
\end{split}
\end{align}
where $u$ denotes unknown complex variable, and the coefficient functions $A(t),B(t)$ and $C(t)$ satisfy the relations:
\begin{equation}\label{relaChazyI}
  \left\{
  \begin{aligned}
  \frac{d^2 A}{dt^2}  &=6A^2,\\
  \frac{d^2 B}{dt^2}  &=6AB,\\
  \frac{d^2 C}{dt^2}  &=B^2+2AC.
   \end{aligned}
  \right. 
\end{equation}
At first, making a change of variables
\begin{equation}
x=u, \quad y=\frac{du}{dt}, \quad z=\frac{d^2u}{dt^2},
\end{equation}
we obtain the system of the first-order ordinary differential equations:
\begin{equation}\label{system;ChazyI}
  \left\{
  \begin{aligned}
   \frac{dx}{dt} &=y,\\
   \frac{dy}{dt} &=z,\\
   \frac{dz}{dt} &=6\left\{-y^2+A(t)\left(y+x^2 \right)+B(t)x+C(t)\right\}.
   \end{aligned}
  \right. 
\end{equation}
Let us make its phase space by gluing two copies of ${\Bbb C}^3 \times {\Bbb C}$ via the birational transformations.

The following Lemma shows that this rational vector field $\tilde v$ associated with the system \eqref{system;ChazyI} has seven accessible singular points on the boundary divisor ${\mathcal H} \subset {\Bbb P}^3$.
\begin{lemma}
The rational vector field $\tilde v$ associated with the system \eqref{system;ChazyI} has three accessible singular points{\rm : \rm}
\begin{equation}
  \left\{
  \begin{aligned}
   P_1 &=\{(X_1,Y_1,Z_1)|X_1=Z_1=0,Y_1=\sqrt{A(t)}\},\\
   P_2 &=\{(X_1,Y_1,Z_1)|X_1=Z_1=0,Y_1=-\sqrt{A(t)}\},\\
   P_3 &=\{(X_3,Y_3,Z_3)|X_3=Y_3=Z_3=0\},
   \end{aligned}
  \right. 
\end{equation}
\end{lemma}
where the point $P_3$ has multiplicity of order 5 and $(X_i,Y_i,Z_i)$ are given by \eqref{P2cover}.

This lemma can be proven by a direct calculation. \qed

Next let us calculate its local index at each point.
\begin{center}
\begin{tabular}{|c|c|c|} \hline 
Singular point & Type of local index   \\ \hline 
$P_1$ & $(0,0,0)$  \\ \hline 
$P_2$ & $(0,0,0)$  \\ \hline
\end{tabular}
\end{center}

\begin{example}
Let us take the coordinate system $(p,q,r)$ centered at the point $P_1$:
$$
p=X_1=\frac{1}{x}, \quad q=Y_1-\sqrt{A(t)}=\frac{y}{x}-\sqrt{A(t)}, \quad r=Z_1=\frac{z}{x}.
$$
The system \eqref{system;ChazyI} is rewritten as follows:
\begin{align*}
\frac{d}{dt}\begin{pmatrix}
             p \\
             q \\
             r 
             \end{pmatrix}&=\frac{1}{p}\left\{\begin{pmatrix}
             0 & 0 & 0 \\
             -A(t)-\frac{A'(t)}{2\sqrt{A(t)}} & 0 & 0 \\
             6(A(t)\sqrt{A(t)}+B(t)) & -12\sqrt{A(t)} & 0
             \end{pmatrix}\begin{pmatrix}
             p \\
             q \\
             r 
             \end{pmatrix}+\cdots\right\}
             \end{align*}
satisfying \eqref{b}. In this case, the local index is $(0,0,0)$. We see that the residue of the formal Laurent series:
\begin{equation}
x(t)=\frac{a_0}{(t-t_0)}+a_1+a_2(t-t_0)+\dots+a_n(t-t_0)^{n-1}+\cdots \quad (a_i \in {\Bbb C})
\end{equation}
is equal to $a_0=0$. By a direct calculation, we see that there are no solutions which pass through $P_1$.
\end{example}
In the case of $P_2$, we also discuss the same way.

In order to do analysis for the accessible singular point $P_3$, we need to replace a suitable coordinate system because this point has multiplicity of order 5.

At first, let us do the Painlev\'e test. To find the leading order behaviour of a singularity at $t=t_0$ one sets
\begin{equation*}
  \left\{
  \begin{aligned}
   x & \propto \frac{a}{(t-t_0)^m},\\
   y & \propto \frac{b}{(t-t_0)^n},\\
   z & \propto \frac{c}{(t-t_0)^p},
   \end{aligned}
  \right. 
\end{equation*}
from which it is easily deduced that
\begin{equation*}
m=1, \quad n=2, \quad p=3.
\end{equation*}
Each order of pole $(m,n,p)$ suggests a suitable coordinate system to do analysis for the accessible singular point $P_3$, which is explicitly given by
\begin{equation*}
(X,Y,Z)=\left(\frac{1}{x},\frac{y}{x^2},\frac{z}{x^3} \right).
\end{equation*}
In this coordinate, the singular points are given as follows:
\begin{equation*}
   P_3^{(1)}= \left\{(X,Y,Z)=\left(0,-1,2 \right)\right\}.
\end{equation*}
Next let us calculate its local index at the point $P_3^{(1)}$.
\begin{center}
\begin{tabular}{|c|c|c|} \hline 
Singular point & Type of local index   \\ \hline 
$P_3^{(1)}$ & $(1,1,6)$  \\ \hline 
\end{tabular}
\end{center}

Now, we try to resolve the accessible singular point $P_3^{(1)}$.

{\bf Step 0}: We take the coordinate system centered at $P_3^{(1)}${\rm : \rm}
$$
p=X, \quad q=Y+1, \quad r=Z-2.
$$

In this coordinate, the system \eqref{system;ChazyI} is rewritten as follows:
\begin{align*}
\frac{d}{dt}\begin{pmatrix}
             p \\
             q \\
             r 
             \end{pmatrix}&=\frac{1}{p}\left\{\begin{pmatrix}
             1 & 0 & 0 \\
             0 & 4 & 1 \\
             0 & 6 & 3
             \end{pmatrix}\begin{pmatrix}
             p \\
             q \\
             r
             \end{pmatrix}+\dots\right\}
             \end{align*}
satisfying \eqref{b}. To the above system, we make the linear transformation
\begin{equation*}
\begin{pmatrix}
             X \\
             Y \\
             Z 
             \end{pmatrix}=\begin{pmatrix}
             1 & 0 & 0 \\
             0 & -\frac{2}{5} & \frac{1}{5} \\
             0 & \frac{3}{5} & \frac{1}{5}
             \end{pmatrix}\begin{pmatrix}
             p \\
             q \\
             r 
             \end{pmatrix}
\end{equation*}
to arrive at
\begin{equation*}
\frac{d}{dt}\begin{pmatrix}
             X \\
             Y \\
             Z 
             \end{pmatrix}=\frac{1}{X}\left\{\begin{pmatrix}
             1 & 0 & 0 \\
             0 & 1 & 0 \\
             0 & 0 & 6 
             \end{pmatrix}\begin{pmatrix}
             X \\
             Y \\
             Z 
             \end{pmatrix}+\cdots\right\}.
             \end{equation*}
By considering the ratio of the local index $(1,1,6)$, we obtain the resonances $\left(\frac{1}{1},\frac{6}{1} \right)=(1,6)$. This property suggests that we will blow up one time to the direction $q$ and six times to the direction $r$.

{\bf Step 1}: We blow up at the point $P_3^{(1)}${\rm : \rm}
$$
p^{(1)}=p, \quad q^{(1)}=\frac{q}{p}, \quad r^{(1)}=\frac{r}{p}.
$$

{\bf Step 2}: We blow up along the curve $\{(p^{(1)},q^{(1)},r^{(1)})=(p^{(1)},q^{(1)},-3q^{(1)})\}${\rm : \rm}
$$
p^{(2)}=p^{(1)}, \quad q^{(2)}=q^{(1)}, \quad r^{(2)}=\frac{r^{(1)}+3q^{(1)}}{p^{(1)}}.
$$

{\bf Step 3}: We blow up along the curve $\left\{(p^{(2)},q^{(2)},r^{(2)})=(p^{(2)},q^{(2)},\frac{3}{4}(q^{(2)})^2)\right\}${\rm : \rm}
$$
p^{(3)}=p^{(2)}, \quad q^{(3)}=q^{(2)}, \quad r^{(3)}=\frac{r^{(2)}-\frac{3}{4}(q^{(2)})^2}{p^{(2)}}.
$$

{\bf Step 4}: We blow up along the curve$\{(p^{(3)},q^{(3)},r^{(3)})=(p^{(3)},q^{(3)},$

$\frac{1}{8}((q^{(3)})^3-16A(t)q^{(3)}-16B(t))\}${\rm : \rm}
$$
p^{(4)}=p^{(3)}, \quad q^{(4)}=q^{(3)}, \quad r^{(4)}=\frac{r^{(3)}-\frac{1}{8}((q^{(3)})^3-16A(t)q^{(3)}-16B(t))}{p^{(3)}}.
$$

{\bf Step 5}: We blow up along the curve $\{(p^{(4)},q^{(4)},r^{(4)})=(p^{(4)},q^{(4)},$

$\frac{1}{64}(3(q^{(4)})^4-80A(t)(q^{(4)})^2-32(2A'(t)+3B(t))q^{(4)}-64B'(t)-192C(t))\}${\rm : \rm}
\begin{align*}
&p^{(5)}=p^{(4)}, \quad q^{(5)}=q^{(4)},\\
&r^{(5)}=\frac{r^{(4)}-\frac{1}{64}(3(q^{(4)})^4-80A(t)(q^{(4)})^2-32(2A'(t)+3B(t))q^{(4)}-64B'(t)-192C(t))}{p^{(4)}}.
\end{align*}

{\bf Step 6}: We blow up along the curve $\{(p^{(5)},q^{(5)},r^{(5)})=(p^{(5)},q^{(5)},r^{(5)}-g(q^{(5)}))${\rm : \rm}
\begin{align*}
&p^{(6)}=p^{(5)}, \quad q^{(6)}=q^{(5)}, \quad r^{(6)}=\frac{r^{(5)}-g(q^{(5)})}{p^{(5)}},
\end{align*}
where the symbol $g(q^{(5)})$ is given by
\begin{align*}
g(q^{(5)})=&\frac{1}{128}(3(q^{(5)})^5-128A(t)(q^{(4)})^3-48(3B(t)+4)(q^{(4)})^2+64(8{A(t)}^2\\
&-3C(t)A'(t)-4B'(t)-2A''(t))q^{(5)}-384C'(t)-128B''(t)+512A(t)B(t)).
\end{align*}
Now, we have blowed up one time to the direction $q$ and six times to the direction $r$. In this coordinate, the system \eqref{system;ChazyI} is rewritten as follows:
\begin{equation}\label{A}
  \left\{
  \begin{aligned}
   \frac{dp^{(6)}}{dt} &=g_1(p^{(6)},q^{(6)},r^{(6)}),\\
   \frac{dq^{(6)}}{dt} &=g_2(p^{(6)},q^{(6)},r^{(6)}),\\
   \frac{dr^{(6)}}{dt} &=\frac{h_1(t)({q^{(6)}})^2+h_2(t)q^{(6)}+h_3(t)}{p^{(6)}}+g_3(p^{(6)},q^{(6)},r^{(6)}),
   \end{aligned}
  \right. 
\end{equation}
where $g_i(p^{(6)},q^{(6)},r^{(6)}) \in {\Bbb C}(t)[p^{(6)},q^{(6)},r^{(6)}] \ (i=1,2,3)$ and $h_1(t)({q^{(6)}})^2+h_2(t)q^{(6)}+h_3(t) \in {\Bbb C}(t)[q^{(6)}]$.

Each right-hand side of the system \eqref{A} is a {\it polynomial} if and only if
\begin{equation}\label{relaChazyI2}
  \left\{
  \begin{aligned}
  \frac{d^2 A}{dt^2}  &=6A^2,\\
  \frac{d^2 B}{dt^2}  &=\frac{2}{3}(9A(t)B(t)+12A'(t)A(t)-A'''(t)),\\
  \frac{d^2 C}{dt^2}  &=\frac{1}{3}(3B(t)^2+6A(t)C(t)+6A'(t)B(t)+6A(t)B'(t)-B'''(t)).
   \end{aligned}
  \right. 
\end{equation}
By solving these equations, we can obtain the relations \eqref{relaChazyI}.

Thus, we have completed the proof of the following theorem.

\begin{theorem}
The phase space ${\mathcal X}$ for the system \eqref{system;ChazyI} is obtained by gluing two copies of ${\Bbb C}^3 \times {\Bbb C}${\rm:\rm}
\begin{center}
${U_j} \times {\Bbb C}={\Bbb C}^3 \times {\Bbb C} \ni \{(x_j,y_j,z_j,t)\},  \ \ j=0,1$
\end{center}
via the following birational transformations{\rm:\rm}
\begin{align}
\begin{split}
0) \ &x_0=x, \quad y_0=y, \quad z_0=z,\\
1) \ &x_1=\frac{1}{x}, \quad y_1=\frac{y+x^2}{x},\\
&z_1=-(((\left(\left(\left(\frac{z}{x^3}-2 \right)x+3\frac{y+x^2}{x} \right)x-\frac{3}{4}\left(\frac{y+x^2}{x}\right)^2 \right)x\\
&-\frac{1}{8}\left(\left(\frac{y+x^2}{x}\right)^3-16A(t)\frac{y+x^2}{x}-16B(t) \right))x\\
&-\frac{1}{64}\left(3\left(\frac{y+x^2}{x}\right)^4-80A(t)\left(\frac{y+x^2}{x}\right)^2-32(2A'(t)+3B(t))\frac{y+x^2}{x}-64B'(t)-192C(t) \right))x\\
&-\frac{1}{128}\{3\left(\frac{y+x^2}{x}\right)^5-128A(t)\left(\frac{y+x^2}{x}\right)^3-48(3B(t)+4)\left(\frac{y+x^2}{x}\right)^2\\
&+64(8{A(t)}^2-3C(t)A'(t)-4B'(t)-2A''(t))\frac{y+x^2}{x}-384C'(t)-128B''(t)+512A(t)B(t)\})x.
\end{split}
\end{align}
\end{theorem}
We remark that these transition functions satisfy the condition{\rm:\rm}
\begin{equation*}
dx_1 \wedge dy_1 \wedge dz_1=dx \wedge dy \wedge dz.
\end{equation*}

\section{Appendix B}
The Chazy-VIII equation is given by
\begin{align}\label{eq;ChazyVIII}
\begin{split}
\frac{d^3 u}{dt^3}=&6u^2 \frac{du}{dt}+(-2\alpha^2 t^2+\beta t+\gamma)\left(\frac{du}{dt}+\alpha \right)+2\alpha u^2+(-4\alpha^2 t+\beta)u,
\end{split}
\end{align}
where $u$ denotes unknown complex variable and $\alpha,\beta$ and $\gamma$ are its constant parameters.

In \cite{Cos1}, we distinguish two canonical subcases. If $\alpha \not=0$, we may set $\beta=0$ by a translation in $t$. This gives equation VIII.a, whose solution is
\begin{equation}
u=w(t)-\alpha t,
\end{equation}
where $w(t)$ satisfies the Painlev\'e IV equation
\begin{equation}
w''=\frac{1}{2w}(w')^2+\frac{3}{2}w^3-4\alpha t w^2+\frac{1}{2}(4\alpha^2 t^2+\gamma)w+\frac{K}{w} \quad \left(':=\frac{d}{dt} \right),
\end{equation}
K being a constant of integration.

If $\alpha=0$, we call the equation Chazy-VIII.b. A first integral of VIII.b is the Painlev\'e II equation:
\begin{equation}
u''=2u^3+(\beta t+\gamma)u+K \quad \left(':=\frac{d}{dt} \right),
\end{equation}
K being a constant of integration.

\begin{theorem}\label{th:appendix}
The birational transformation $\varphi_0$
\begin{equation}\label{bira:A0}
  \left\{
  \begin{aligned}
   x &=u,\\
   y &=\frac{du}{dt}-u^2-\frac{1}{2}\{(-2\alpha^2 t+\beta)t+2\alpha+\gamma \},\\
   z &=\frac{d^2u}{dt^2}-2u^3-\{(-2\alpha^2 t+\beta)t+2\alpha+\gamma \} u
   \end{aligned}
  \right. 
\end{equation}
takes the Chazy VIII equation to the system of the first-order ordinary differential equations:
\begin{equation}\label{system;ChazyVIII}
  \left\{
  \begin{aligned}
   \frac{dx}{dt} &=x^2+y+\frac{\beta}{2}t-\alpha^2 t^2+\alpha+\frac{\gamma}{2},\\
   \frac{dy}{dt} &=-2xy+z+2\alpha^2 t-\frac{\beta}{2},\\
   \frac{dz}{dt} &=-2\alpha y-2\alpha^2.
   \end{aligned}
  \right. 
\end{equation}
\end{theorem}

By the same way of the second Painlev\'e system, we can prove Theorem \ref{th:appendix}.

{\bf Proof.} At first, we rewrite the equation \eqref{eq;ChazyVIII} to the system of the first-order ordinary differential equations.

{\bf Step 0:} We make a change of variables.
\begin{equation}
x=u, \quad y=\frac{du}{dt}, \quad z=\frac{d^2u}{dt^2}.
\end{equation}

{\bf Step 1:} We make a change of variables.
\begin{equation}
x_1=\frac{1}{x}, \quad y_1=\frac{y}{x^2}, \quad z_1=\frac{z}{x^3}.
\end{equation}
In this coordinate system, we see that this system has two accessible singular points:
\begin{equation}
(x_1,y_1,z_1)=\left\{(0,1,2),(0,-1,2) \right\}.
\end{equation}

Around the point $(x_1,y_1,z_1)=(0,1,2)$, we can rewrite the system as follows.

{\bf Step 2:} We make a change of variables.
\begin{equation}
x_2=x_1, \quad y_2=y_1-1, \quad z_2=z_1-2.
\end{equation}
In this coordinate system, we can rewrite the system satisfying the condition \eqref{b}:
\begin{align*}
\frac{d}{dt}\begin{pmatrix}
             x_2 \\
             y_2 \\
             z_2 
             \end{pmatrix}&=\frac{1}{x_2}\left\{\begin{pmatrix}
             -1 & 0 & 0  \\
             0 & -4 & 1  \\
             0 & 0 & -3
             \end{pmatrix}\begin{pmatrix}
             x_2 \\
             y_2 \\
             z_2 
             \end{pmatrix}+\cdots\right\}.
             \end{align*}
and we can obtain the local index $(-1,-4,-3)$ at the point $\{(x_2,y_2,z_2)=(0,0,0)\}$. The continued ratio of the local index at the point $\{(x_2,y_2,z_2)=(0,0,0)\}$ are all positive integers
\begin{equation}
\left(\frac{-4}{-1},\frac{-3}{-1} \right)=(4,3).
\end{equation}
This is the reason why we choose this accessible singular point.

We aim to obtain the local index $(-1,-2,0)$ by successive blowing-up procedures.

{\bf Step 3:} We blow up at the point $\{(x_2,y_2,z_2)=(0,0,0)\}$.
\begin{equation}
x_3=x_2, \quad y_3=\frac{y_2}{x_2}, \quad z_3=\frac{z_2}{x_2}.
\end{equation}

{\bf Step 4:} We blow up at the point $\{(x_3,y_3,z_3)=(0,0,0)\}$.
\begin{equation}
x_4=x_3, \quad y_4=\frac{y_3}{x_3}, \quad z_4=\frac{z_3}{x_3}.
\end{equation}
In this coordinate system, we see that this system has the following accessible singular point:
\begin{equation}
(x_4,y_4,z_4)=\left(0,\frac{1}{2}(2\alpha-2\alpha^2 t^2+\beta t+\gamma),2\alpha-2\alpha^2 t^2+\beta t+\gamma \right).
\end{equation}

{\bf Step 5:} We blow up along the curve $\{(x_4,y_4,z_4)=(0,\frac{1}{2}(2\alpha-2\alpha^2 t^2+\beta t+\gamma),\\
2\alpha-2\alpha^2 t^2+\beta t+\gamma )\}$.
\begin{equation}
x_5=x_4, \quad y_5=y_4-\frac{1}{2}(2\alpha-2\alpha^2 t^2+\beta t+\gamma), \quad z_5=\frac{z_4-(2\alpha-2\alpha^2 t^2+\beta t+\gamma)}{x_4}.
\end{equation}

{\bf Step 6:} We make a change of variables.
\begin{equation}
x_6=\frac{1}{x_5}, \quad y_6=y_5, \quad z_6=z_5,
\end{equation}
and we can obtain the system \eqref{system;ChazyVIII}.

Thus, we have completed the proof of Theorem \ref{th:appendix}. \qed

We note on the remaining accessible singular points.

Around the point $(x_1,y_1,z_1)=(0,-1,2)$, we can rewrite the system as follows.

{\bf Step 2:} We make a change of variables.
\begin{equation}
x_2=x_1, \quad y_2=y_1+1, \quad z_2=z_1-2.
\end{equation}
\begin{align*}
\frac{d}{dt}\begin{pmatrix}
             x_2 \\
             y_2 \\
             z_2 
             \end{pmatrix}&=\frac{1}{x_2}\left\{\begin{pmatrix}
             1 & 0 & 0  \\
             0 & 4 & 1  \\
             0 & 0 & 3
             \end{pmatrix}\begin{pmatrix}
             x_2 \\
             y_2 \\
             z_2 
             \end{pmatrix}+\cdots\right\}.
             \end{align*}
and we can obtain the local index $(1,4,3)$ at the point $\{(x_2,y_2,z_2)=(0,0,0)\}$. The continued ratio of the local index at the point $\{(x_2,y_2,z_2)=(0,0,0)\}$ are all positive integers
\begin{equation}
\left(\frac{4}{1},\frac{3}{1} \right)=(4,3).
\end{equation}
We remark that we can discuss this case in the same way as in the case of $(x_1,y_1,z_1)=(0,1,2)$.

Now, let us consider the birational B{\"a}cklund transformations of the Chazy-VIII system. These B{\"a}cklund transformations are new.

\begin{proposition}
The birational transformation $s_0$
\begin{equation}\label{bira:A1}
  \left\{
  \begin{aligned}
   X &=x-\frac{2z+4\alpha^2 t-\beta}{2y},\\
   Y &=y,\\
   Z &=-z
   \end{aligned}
  \right. 
\end{equation}
takes the system \eqref{system;ChazyVIII} to the system
\begin{equation}\label{system;A1}
  \left\{
  \begin{aligned}
   \frac{dX}{dt} &=X^2+Y-\alpha^2 t^2+\frac{\beta}{2}t+3\alpha+\frac{\gamma}{2},\\
   \frac{dY}{dt} &=Z-2XY-2\alpha^2 t+\frac{\beta}{2},\\
   \frac{dZ}{dt} &=2\alpha Y+2\alpha^2.
   \end{aligned}
  \right. 
\end{equation}
\end{proposition}

\begin{proposition}
The birational transformation $s_1$
\begin{equation}\label{bira:A2}
  \left\{
  \begin{aligned}
   X &=-\left(x+\frac{z+4\alpha x-2\alpha^2 t+\frac{\beta}{2}}{y+2x^2-2\alpha^2 t^2+\beta t+\gamma} \right),\\
   Y &=-(y+2x^2-2\alpha^2 t^2+\beta t+\gamma),\\
   Z &=-(-z-4\alpha x)
   \end{aligned}
  \right. 
\end{equation}
takes the system \eqref{system;ChazyVIII} to the system
\begin{equation}\label{system;A2}
  \left\{
  \begin{aligned}
   \frac{dX}{dt} &=X^2+Y-\alpha^2 t^2+\frac{\beta}{2}t-3\alpha+\frac{\gamma}{2},\\
   \frac{dY}{dt} &=Z-2XY-2\alpha^2 t+\frac{\beta}{2},\\
   \frac{dZ}{dt} &=-2\alpha Y+2\alpha^2.
   \end{aligned}
  \right. 
\end{equation}
\end{proposition}
Here, its inverse transformation ${s_1}^{-1}$ is given by
\begin{equation}\label{bira:A3}
  \left\{
  \begin{aligned}
   x =&-X+\frac{2Z-4\alpha^2 t+\beta}{2Y},\\
   y =&-Y-2X^2+2\alpha^2 t^2-\beta t-\gamma+\frac{4X(Z-2\alpha^2 t)}{Y}\\
   &+\frac{4\beta XY-4Z^2+4Z(4\alpha^2 t-\beta)-16\alpha^4 t^2+8\alpha^2 \beta t-\beta^2}{2Y^2},\\
   z =&Z+4\alpha X-\frac{2\alpha(2Z-4\alpha^2 t+\beta)}{Y}.
   \end{aligned}
  \right. 
\end{equation}

\begin{proposition}
The transformation
\begin{equation}\label{bira:A4}
\varphi:(X,Y,Z;\alpha,\beta,\gamma) \rightarrow (X,Y,Z;-\alpha,\beta,\gamma)
\end{equation}
is a B{\"a}cklund transformation from the system \eqref{system;A1} to the system \eqref{system;A2}.
\end{proposition}

\begin{figure}[t]
\unitlength 0.1in
\begin{picture}( 31.7500, 12.4000)( 16.9000,-19.3000)
\put(16.9000,-9.9000){\makebox(0,0)[lb]{$(x,y,z)$}}%
%
\special{pn 20}%
\special{pa 2066 1090}%
\special{pa 2066 1760}%
\special{fp}%
\special{sh 1}%
\special{pa 2066 1760}%
\special{pa 2086 1694}%
\special{pa 2066 1708}%
\special{pa 2046 1694}%
\special{pa 2066 1760}%
\special{fp}%
\put(16.9000,-19.7000){\makebox(0,0)[lb]{$(X,Y,Z)$}}%
%
\special{pn 20}%
\special{pa 2630 1870}%
\special{pa 4260 1870}%
\special{fp}%
\special{sh 1}%
\special{pa 4260 1870}%
\special{pa 4194 1850}%
\special{pa 4208 1870}%
\special{pa 4194 1890}%
\special{pa 4260 1870}%
\special{fp}%
\put(44.9000,-10.1000){\makebox(0,0)[lb]{$(x,y,z)$}}%
\put(44.9000,-19.9000){\makebox(0,0)[lb]{$(X,Y,Z)$}}%
%
\special{pn 20}%
\special{pa 4826 1760}%
\special{pa 4826 1100}%
\special{fp}%
\special{sh 1}%
\special{pa 4826 1100}%
\special{pa 4806 1168}%
\special{pa 4826 1154}%
\special{pa 4846 1168}%
\special{pa 4826 1100}%
\special{fp}%
%
\special{pn 8}%
\special{pa 2630 900}%
\special{pa 4240 900}%
\special{fp}%
\special{sh 1}%
\special{pa 4240 900}%
\special{pa 4174 880}%
\special{pa 4188 900}%
\special{pa 4174 920}%
\special{pa 4240 900}%
\special{fp}%
\put(31.7000,-8.6000){\makebox(0,0)[lb]{$\pi$}}%
\put(20.9500,-15.1000){\makebox(0,0)[lb]{$s_0$}}%
\put(48.6500,-15.1000){\makebox(0,0)[lb]{${s_1}^{-1}$}}%
\put(31.9000,-21.0000){\makebox(0,0)[lb]{$\varphi$}}%
\end{picture}%
\label{ChazyVIIIfig1}
\caption{The transformation $\pi$ can be obtained by the compositions of the transformations \eqref{bira:A1}, \eqref{bira:A4} and  \eqref{bira:A3}.}
\end{figure}

\begin{theorem}\label{prp:A4.4}
The system \eqref{system;ChazyVIII} is invariant under the following birational transformations$:$
\begin{align}
\begin{split}
&s_0:(x,y,z;0,0,\gamma) \rightarrow \left(x-\frac{z}{y},y,-z;0,0,\gamma \right),\\
&s_1:(x,y,z;0,0,\gamma) \rightarrow \left\{-\left(x+\frac{z}{y+2x^2+\gamma} \right),-(y+2x^2+\gamma),z;0,0,\gamma \right\},\\
&\pi:(x,y,z;\alpha,\beta,\gamma) \rightarrow (-x,-y-2x^2+2\alpha^2 t^2-\beta t-\gamma,-z-4\alpha x;-\alpha,\beta,\gamma).
\end{split}
\end{align}
\end{theorem}
These transformations are new. It should be clear from the form of $s_0$ and $s_1$ that the transformations become auto-B{\"a}cklund transformations for the system \eqref{system;ChazyVIII} only if $\alpha=\beta=0$. The transformation $\pi$ can be obtained by the compositions of the transformations \eqref{bira:A1}, \eqref{bira:A4} and  \eqref{bira:A3}.

Next, we give the holomorphy condition of the system \eqref{system;ChazyVIII}. Thanks to this holomorphy condition, we can recover the system \eqref{system;ChazyVIII}.
\begin{theorem}
Let us consider the following ordinary differential system in the polynomial class\rm{:\rm}
\begin{equation*}
  \left\{
  \begin{aligned}
   \frac{dx}{dt} &=f_1(x,y,z),\\
   \frac{dy}{dt} &=f_2(x,y,z),\\
   \frac{dz}{dt} &=f_3(x,y,z).
   \end{aligned}
  \right. 
\end{equation*}
We assume that

$(A1)$ $deg(f_i)=3$ with respect to $x,y,z$.

$(A2)$ The right-hand side of this system becomes again a polynomial in each coordinate system $(x_i,y_i,z_i) \ (i=1,2)$.
\begin{align}
\begin{split}
1) \ &x_1=\frac{1}{x}, \quad y_1=-\left(yx-\frac{1}{2}(2z+4\alpha^2 t-\beta) \right)x, \quad z_1=z,\\
2) \ &x_2=\frac{1}{x}, \quad y_2=-\left((y+2x^2-2\alpha^2 t^2+\beta t+4\alpha+\gamma)x+z-2\alpha^2 t+\frac{\beta}{2} \right)x,\\
&z_2=z+4\alpha x.
\end{split}
\end{align}
Then such a system coincides with the system \eqref{system;ChazyVIII}.
\end{theorem}
These transition functions satisfy the condition{\rm:\rm}
\begin{equation*}
dx_i \wedge dy_i \wedge dz_i=dx \wedge dy \wedge dz \quad (i=1,2).
\end{equation*}

Finally, we study a solution of the system \eqref{system;ChazyVIII} which is written by the use of known functions.

\begin{proposition}
The system \eqref{system;ChazyVIII} admits a rational solution{\rm:\rm}
\begin{align}
\begin{split}
&(x,y,z;\alpha,\beta,\gamma)=\left(0,-\frac{\beta}{2}t-\frac{\gamma}{2},0;0,\beta,\gamma \right),
\end{split}
\end{align}
and  admits a special solution{\rm:\rm}
\begin{align}\label{rat}
\begin{split}
(x,y,z;\alpha,\beta,\gamma)=&\left(x,0,\frac{1}{2}(\beta-4\alpha^2 t);\alpha,\beta,\gamma \right),
\end{split}
\end{align}
where the equation in $x$ satisfies as follows:
\begin{equation}
\frac{dx}{dt}=x^2+\frac{\beta}{2}t-\alpha^2 t^2+\alpha+\frac{\gamma}{2}.
\end{equation}
\end{proposition}
Setting $x:=-\frac{d}{dt}log X=-\frac{X'}{X}$, this equation can be transformed into the second-order linear differential equation:
\begin{equation}
\frac{d^2X}{dt^2}=\frac{1}{2}(2\alpha^2 t^2-\beta t-2\alpha-\gamma)X.
\end{equation}
This equation can be solved by the ParabolicCylinder functions.

We remark that the solution \eqref{rat} also can be obtained as the fixed point of the B{\"a}cklund transformation $\pi$ (see Proposition \ref{prp:A4.4}).

\section{Appendix C}
In this section, we present a 6-parameter family of ordinary differential systems in dimension three explicitly given by
\begin{equation}\label{mmSVIII}
  \left\{
  \begin{aligned}
   \frac{dx}{dt} &=x^2-xy-xz+(-\alpha_3+\alpha_4-\alpha_5+\alpha_6)x+\alpha_3 y+\alpha_5 z+\alpha_3 \alpha_5-\alpha_4 \alpha_5 -\alpha_3 \alpha_6,\\
   \frac{dy}{dt} &=y^2-xy-yz+\alpha_1 x+(-\alpha_1+\alpha_2+\alpha_5-\alpha_6)y+\alpha_6 z-\alpha_1 \alpha_5+\alpha_1 \alpha_6-\alpha_2 \alpha_6,\\
   \frac{dz}{dt} &=z^2-xz-yz+\alpha_2 x+\alpha_4 y+(\alpha_1-\alpha_2+\alpha_3-\alpha_4)z-\alpha_2 \alpha_3-\alpha_1 \alpha_4+\alpha_2 \alpha_4.
   \end{aligned}
  \right. 
\end{equation}
Here $x,y,z$ denote unknown complex variables and $\alpha_i \ (i=1,2,\ldots,6)$ are complex parameters.

\begin{proposition}
This system is invariant under the following transformations{\rm:\rm}
\begin{align}
\begin{split}
s_0(x,y,z;\alpha_1,\ldots,\alpha_6) \rightarrow &(y,x,z;\alpha_3,\alpha_4,\alpha_1,\alpha_2,\alpha_6,\alpha_5),\\
s_1(x,y,z;\alpha_1,\ldots,\alpha_6) \rightarrow &(z,y,x;\alpha_6,\alpha_5,\alpha_4,\alpha_3,\alpha_2,\alpha_1),\\
s_2(x,y,z;\alpha_1,\ldots,\alpha_6) \rightarrow &(x,z,y;\alpha_2,\alpha_1,\alpha_5,\alpha_6,\alpha_3,\alpha_4),\\
\pi(x,y,z;\alpha_1,\ldots,\alpha_6) \rightarrow &(y,z,x;\alpha_4,\alpha_3,\alpha_6,\alpha_5,\alpha_1,\alpha_2).
\end{split}
\end{align}
\end{proposition}

\begin{theorem}
After a series of explicit blowing-ups at eight points including four infinitely near points on the boundary divisor ${\mathcal H} \cong {\Bbb P}^2$ in ${\Bbb P}^3$, we obtain the smooth projective 3-fold $\tilde{\mathcal X}$ and a morphism $\varphi:\tilde{\mathcal X} \rightarrow {\Bbb P}^3$. Its canonical divisor $K_{\tilde{\mathcal X}}$ is given by
\begin{align}
\begin{split}
K_{\tilde{\mathcal X}}&=-4{\mathcal E}_0 -2\sum_{i=1}^{4} {\mathcal E}_i,
\end{split}
\end{align}
where the symbol ${\mathcal E}_0$ denotes the proper transform of boundary divisor ${\mathcal H}$ of ${\Bbb P}^3$ by $\varphi$ and  ${\mathcal E}_i$ denote the exceptional divisors, which are isomorphic to ${\Bbb F}_1$. Moreover, $\tilde{\mathcal X}-(-{K_{\tilde{\mathcal X}}})_{red}$ satisfies
\begin{equation}
\tilde{\mathcal X}-(-{K_{\tilde{\mathcal X}}})_{red}={\mathcal X}.
\end{equation}
\end{theorem}

\begin{theorem}
The phase space ${\mathcal X}$ for the system \eqref{mmSVIII} is obtained by gluing five copies of ${\Bbb C}^3${\rm:\rm}
\begin{center}
${U_j} \cong {\Bbb C}^3 \ni \{(x_j,y_j,z_j)\},  \ \ j=0,1,\ldots,4$
\end{center}
via the following birational transformations{\rm:\rm}
\begin{align}
\begin{split}
0) \ &x_0=x, \quad y_0=y, \quad z_0=z,\\
1) \ &x_1=\frac{1}{x}, \quad y_1=-(y-\alpha_1)x, \quad z_1=(z-\alpha_2)x,\\
2) \ &x_2=(x-\alpha_3)y, \quad y_2=\frac{1}{y}, \quad z_2=-(z-\alpha_4)y,\\
3) \ &x_3=-(x-\alpha_5)z, \quad y_3=(y-\alpha_6)z, \quad z_3=\frac{1}{z},\\
4) \ &x_4=\frac{1}{x}, \quad y_4=-(y-x+\alpha_2-\alpha_4+\alpha_5-\alpha_6)x,\\
&z_4=(z-x+\alpha_1+\alpha_3-\alpha_4-\alpha_6)x.
\end{split}
\end{align}
These transition functions satisfy the condition{\rm:\rm}
\begin{equation*}
dx_i \wedge dy_i \wedge dz_i=dx \wedge dy \wedge dz \quad (i=1,2,3,4).
\end{equation*}
\end{theorem}

\begin{proposition}
This system has
\begin{equation}\label{imm}
I:=xz-yz-\alpha_2 x+\alpha_4 y-(\alpha_5-\alpha_6)z
\end{equation}
as its first integral.
\end{proposition}

By using \eqref{imm}, elimination of $x$ from the system \eqref{mmSVIII} gives the second-order ordinary differential system in the variables $(y,z)$;
\begin{equation}\label{mm-s}
  \left\{
  \begin{aligned}
   \frac{dy}{dt} =&y^2-yz+(-\alpha_1+\alpha_2+\alpha_5-\alpha_6)y+\alpha_6 z-\alpha_1 \alpha_5+\alpha_1 \alpha_6 -\alpha_2\alpha_6\\
   &-\frac{(y-\alpha_1)(I+yz-\alpha_4 y+(\alpha_5-\alpha_6) z)}{z-\alpha_2},\\
   \frac{dz}{dt} =&z^2-2yz+2\alpha_4 y+(\alpha_1-\alpha_2+\alpha_3-\alpha_4-\alpha_5+\alpha_6)z\\
   &-\alpha_1 \alpha_4+\alpha_2 \alpha_4-\alpha_2 \alpha_3-I.
   \end{aligned}
  \right. 
\end{equation}

\begin{proposition}
The canonical transformation
\begin{equation}
  \left\{
  \begin{aligned}
   X &=\frac{y-\alpha_1}{z-\alpha_2},\\
   Y &=z-\alpha_2
   \end{aligned}
  \right. 
\end{equation}
takes the system \eqref{mm-s} to the Hamiltonian system
\begin{equation}
  \left\{
  \begin{aligned}
   \frac{dX}{dt}=\frac{\partial H}{\partial Y}=&2X^2Y+(\alpha_2-\alpha_4)X^2-2XY+(\alpha_1-\alpha_2-\alpha_3+\alpha_4+\alpha_5-\alpha_6)X\\
   &+\alpha_6-\alpha_1,\\
   \frac{dY}{dt}=-\frac{\partial H}{\partial X}=&-2XY^2+Y^2-2(\alpha_2-\alpha_4)XY-(\alpha_1-\alpha_2-\alpha_3+\alpha_4+\alpha_5-\alpha_6)Y\\
   &-I-\alpha_1 \alpha_2+\alpha_1 \alpha_4-\alpha_2 \alpha_5+\alpha_2 \alpha_6
   \end{aligned}
  \right. 
\end{equation}
with the polynomial Hamiltonian
\begin{align*}
H:=&X^2Y^2+(\alpha_2-\alpha_4)X^2Y-XY^2+(\alpha_1-\alpha_2-\alpha_3+\alpha_4+\alpha_5-\alpha_6)XY\\
&-(-I-\alpha_1 \alpha_2+\alpha_1 \alpha_4-\alpha_2 \alpha_5+\alpha_2 \alpha_6)X+(\alpha_6-\alpha_1)Y.
\end{align*}
\end{proposition}
This system is an autonomous version of the fifth Painlev\'e system.

\section{Appendix D}
In this section, we present a 3-parameter family of ordinary differential systems in dimension three  explicitly given by
\begin{equation}\label{system;m-Chazy-VII}
  \left\{
  \begin{aligned}
   \frac{dx}{dt} &=x^2-xy-(\alpha_1-2\alpha_3)x+(\alpha_1-\alpha_3)y-(\alpha_1-\alpha_3)\alpha_3,\\
   \frac{dy}{dt} &=y^2-xy+xz-yz+(\alpha_1-\alpha_2)x-(\alpha_1-\alpha_2+\alpha_3)y+\alpha_3 z+(\alpha_1-\alpha_2)\alpha_3,\\
   \frac{dz}{dt} &=z^2-3xz+3\alpha_2 x+(3\alpha_1-2\alpha_2-3\alpha_3)z-\alpha_2(3\alpha_1-\alpha_2-3\alpha_3).
   \end{aligned}
  \right. 
\end{equation}
Here $x,y,z$ denote unknown complex variables and $\alpha_i \ (i=1,2,3)$ are complex parameters.

\begin{proposition}
This system has
\begin{align}
\begin{split}
I=&2x^3(z-\alpha_2)+x^2\{y^2-2y(z+\alpha_1-\alpha_2)-2(2\alpha_1-3\alpha_3)z-6\alpha_2\alpha_3+\alpha_1^2+4\alpha_1\alpha_2\}\\
&-2x(\alpha_1-\alpha_3)\{y^2-2y(z+\alpha_1-\alpha_2)-(\alpha_1-3\alpha_3)z+\alpha_1^2+\alpha_1\alpha_2-3\alpha_2\alpha_3\}\\
&+(\alpha_1-\alpha_3)^2\{y^2-2y(z+\alpha_1-\alpha_2)+2\alpha_3 z\}
\end{split}
\end{align}
as its first integral.
\end{proposition}

The following Lemma shows that the rational vector field $\tilde v$ associated with the system \eqref{system;m-Chazy-VII} has seven accessible singular points on the boundary divisor ${\mathcal H} \subset {\Bbb P}^3$.
\begin{lemma}
The rational vector field $\tilde v$ associated with the system \eqref{system;m-Chazy-VII} has seven accessible singular points{\rm : \rm}
\begin{equation}
  \left\{
  \begin{aligned}
   P_1 &=\{(X_1,Y_1,Z_1)|X_1=Y_1=Z_1=0\},\\
   P_2 &=\{(X_2,Y_2,Z_2)|X_2=Y_2=Z_2=0\},\\
   P_3 &=\{(X_3,Y_3,Z_3)|X_3=Y_3=Z_3=0\},\\
   P_4 &=\left\{(X_1,Y_1,Z_1)|X_1=0, \ Y_1=\frac{4}{3}, \ Z_1=\frac{8}{3}\right\},\\
   P_5 &=\left\{(X_2,Y_2,Z_2)|X_2=Y_2=0, \ Z_2=\frac{1}{2}\right\},\\
   P_6 &=\{(X_1,Y_1,Z_1)|X_1=0, \ Y_1=1, \ Z_1=3\},\\
   P_7 &=\{(X_1,Y_1,Z_1)|Y_1=1, \ X_1=Z_1=0\},
   \end{aligned}
  \right. 
\end{equation}
where $(X_i,Y_i,Z_i)$ are given by \eqref{P2cover}.
\end{lemma}
This lemma can be proven by a direct calculation. \qed

Next let us calculate its local index at each point.
\begin{center}
\begin{tabular}{|c|c|c|} \hline 
Singular point & Type of local index   \\ \hline 
$P_1$ & $(1,2,4)$  \\ \hline 
$P_2$ & $(2,1,1)$  \\ \hline 
$P_3$ & $(1,2,1)$  \\ \hline 
$P_4$ & $(1,6,4)$  \\ \hline 
$P_5$ & $(3,1,-2)$  \\ \hline 
$P_6$ & $(0,2,-3)$  \\ \hline 
$P_7$ & $(0,-1,2)$  \\ \hline 
\end{tabular}
\end{center}

By resolving the accessible singular points, we can obtain the phase space ${\mathcal X}$ for the system \eqref{system;m-Chazy-VII}.
\begin{theorem}
The phase space ${\mathcal X}$ for the system \eqref{system;m-Chazy-VII} is obtained by gluing six copies of ${\Bbb C}^3${\rm:\rm}
\begin{center}
${U_j}={\Bbb C}^3 \ni \{(x_j,y_j,z_j)\},  \ \ j=0,1,\ldots,5$
\end{center}
via the following birational transformations{\rm:\rm}
\begin{align}
\begin{split}
0) \ &x_0=x, \quad y_0=y, \quad z_0=z,\\
1) \ &x_1=\frac{1}{x}, \quad y_1=(y-\alpha_1)x, \quad z_1=x^3(z-\alpha_2),\\
2) \ &x_2=(x-(\alpha_1-\alpha_3))y, \quad y_2=\frac{1}{y}, \quad z_2=z,\\
3) \ &x_3=x, \quad y_3=(y-x-\alpha_3)z, \quad z_3=\frac{1}{z},\\
4) \ &x_4=\frac{1}{x}, \quad y_4=-\frac{2}{3}x^5(16x-6y-3z-10\alpha_1+3\alpha_2+16\alpha_3),\\
&z_4=x^3(8x-4y-z-4\alpha_1+\alpha_2+8\alpha_3),\\
5) \ &x_5=-(x-y+\alpha_3)(x-\alpha_1+\alpha_3)z, \quad y_5=-\frac{1}{(x-y+\alpha_3)z}, \quad z_5=\frac{1}{z}.
\end{split}
\end{align}
\end{theorem}

\begin{theorem}
Let us consider a system of the first-order ordinary differential equations in the polynomial class\rm{:\rm}
\begin{equation*}
\frac{dx}{dt}=f_1(x,y,z), \quad \frac{dy}{dt}=f_2(x,y,z), \quad \frac{dz}{dt}=f_3(x,y,z).
\end{equation*}
We assume that

$(A1)$ $deg(f_i)=2$ with respect to $x,y,z$.

$(A2)$ The right-hand side of this system becomes again a polynomial in each coordinate system $(x_i,y_i,z_i) \ (i=1,2,3,4)$.

\noindent
Then such a system coincides with the system \eqref{system;m-Chazy-VII}.
\end{theorem}

\section{Appendix E}
The Chazy-XI equation with $N=3$ is explicitly given by
\begin{equation}\label{eq:Chazy-XI}
\frac{d^3 u}{dt^3}=3u^4+6u^2\frac{du}{dt}+\left(\frac{du}{dt}\right)^2-2u\frac{d^2u}{dt^2}.
\end{equation}
Here $u$ denotes unknown complex variable.

\begin{proposition}
The equation \eqref{eq:Chazy-XI} is equivalent to the system of the first-order ordinary differential equations$:$
\begin{equation}\label{system:Chazy-XI}
  \left\{
  \begin{aligned}
   \frac{dx}{dt} &=x^2-2xy-2yz,\\
   \frac{dy}{dt} &=y^2-2xy,\\
   \frac{dz}{dt} &=xz.
   \end{aligned}
  \right. 
\end{equation}
Here $x,y,z$ denote unknown complex variables.
\end{proposition}

\begin{proposition}
This system admits rational solutions\rm{:\rm}
\begin{equation}
  \left\{
  \begin{aligned}
   x(t) &=-\frac{1}{t+c_1},\\
   y(t) &=\frac{3(t^2+2c_1t+c_1^2)}{t^3+3c_1t^2+3c_1^2t+3c_2},\\
   z(t) &=\frac{1}{t+c_1} \quad (c_1,c_2 \in {\Bbb C}),
   \end{aligned}
  \right. 
\end{equation}
and
\begin{equation}
  \left\{
  \begin{aligned}
   x(t) &=-\frac{c_2}{c_2t-c_1},\\
   y(t) &=0,\\
   z(t) &=-\frac{1}{c_2t-c_1} \quad (c_1,c_2 \in {\Bbb C}).
   \end{aligned}
  \right. 
\end{equation}
\end{proposition}

\begin{proposition}
This system has
\begin{equation}
   I=(x+z)y^2z^3
\end{equation}
as its first integral.
\end{proposition}
By using this, elimination of $x$ from the system \eqref{system:Chazy-XI} gives the second-order ordinary differential system for $(y,z)$; namely,
\begin{equation}
  \left\{
  \begin{aligned}
   \frac{dy}{dt} &=y^2-\frac{2(I-y^2z^4)}{yz^3},\\
   \frac{dz}{dt} &=\frac{I-y^2z^4}{y^2z^2}.
   \end{aligned}
  \right. 
\end{equation}
By making a change of the variables
\begin{equation}
  \left\{
  \begin{aligned}
   X &=-\frac{1}{yz^2},\\
   Y &=\frac{1}{z},
   \end{aligned}
  \right. 
\end{equation}
we obtain an autonomous version of Painlev\'e IV system:
\begin{equation}\label{mautoPIV}
  \left\{
  \begin{aligned}
   \frac{dX}{dt} &=Y^2,\\
   \frac{dY}{dt} &=-IX^2+1.
   \end{aligned}
  \right. 
\end{equation}
Now, we give a generalization of the system \eqref{mautoPIV} in addition to constant complex parameters:
\begin{equation}\label{ggggg}
  \left\{
  \begin{aligned}
   \frac{ds}{dt} &=c^2-as+\alpha_1,\\
   \frac{dc}{dt} &=-s^2+ac+\alpha_2.
   \end{aligned}
  \right. 
\end{equation}
Here $s,c$ denote unknown complex variables, and $\alpha_i,a$ are constant complex parameters.

\begin{theorem}
The system \eqref{ggggg} has extended affine Weyl group symmetry of type $A_2^{(1)}$, whose generators $w_i,\pi$ are given by
\begin{align*}\label{D4}
\begin{split}
w_0(*) \rightarrow & \left(s-\frac{\alpha_1+\alpha_2}{s+c+a},c+\frac{\alpha_1+\alpha_2}{s+c+a};-\alpha_2,-\alpha_1 \right),\\
w_1(*) \rightarrow & \left(s-\frac{(-1)^{\frac{1}{3}}((-1)^{\frac{1}{3}} \alpha_1-\alpha_2)}{c-(-1)^{\frac{1}{3}}s+(-1)^{\frac{2}{3}}a},c+\frac{\alpha_1+(-1)^{\frac{2}{3}} \alpha_2}{c-(-1)^{\frac{1}{3}}s+(-1)^{\frac{2}{3}}a};-(-1)^{\frac{2}{3}} \alpha_2,(-1)^{\frac{1}{3}} \alpha_1 \right),\\
w_2(*) \rightarrow & \left(s-\frac{(-1)^{\frac{1}{3}}((-1)^{\frac{1}{3}} \alpha_2-\alpha_1)}{c+(-1)^{\frac{2}{3}}s-(-1)^{\frac{1}{3}}a},c-\frac{(-1)^{\frac{1}{3}} \alpha_2-\alpha_1}{c+(-1)^{\frac{2}{3}}s-(-1)^{\frac{1}{3}}a};(-1)^{\frac{1}{3}} \alpha_2,-(-1)^{\frac{2}{3}} \alpha_1 \right),\\
\pi(*) \rightarrow & \left(-\frac{s}{(-1)^{\frac{1}{3}}},-(-1)^{\frac{1}{3}} c;-(-1)^{\frac{2}{3}} \alpha_1,-(-1)^{\frac{1}{3}} \alpha_2 \right).
\end{split}
\end{align*}
Here the symbol $(*)$ denotes $(s,c;\alpha_1,\alpha_2)$.
\end{theorem}

\begin{theorem}
The phase space $S$ for the system \eqref{ggggg} is a rational surface of type $E_6^{(1)}$, which is obtained by gluing four copies of ${\Bbb C}^2${\rm:\rm}
\begin{center}
$(x_j,y_j) \in {U_j} \cong {\Bbb C}^2 \ \ (j=0,1,2,3),$
\end{center}
via the following birational and symplectic transformations{\rm:\rm}
\begin{align*}
\begin{split}
0) \ &x_0=s, \quad y_0=c,\\
1) \ &x_1=\frac{1}{s}, \quad y_1=-((c+s+a)s-\alpha_1-\alpha_2)s,\\
2) \ &x_2=\frac{1}{s}, \quad y_2=-\{(c-(-1)^{\frac{1}{3}}s+(-1)^{\frac{2}{3}}a)s-(-1)^{\frac{1}{3}}((-1)^{\frac{1}{3}} \alpha_1-\alpha_2)\}s,\\
3) \ &x_3=\frac{1}{s}, \quad y_3=-\{(c+(-1)^{\frac{2}{3}}s-(-1)^{\frac{1}{3}}a)s-(-1)^{\frac{2}{3}}((-1)^{\frac{2}{3}} \alpha_1+\alpha_2)\}s.
\end{split}
\end{align*}
\end{theorem}

The following Lemma shows that this rational vector field $\tilde v$ associated with the system \eqref{system:Chazy-XI} has four accessible singular loci on the boundary divisor ${\mathcal H} \subset {\Bbb P}^3$.
\begin{lemma}
The rational vector field $\tilde v$ associated with the system \eqref{system:Chazy-XI} has four accessible singular loci{\rm : \rm}
\begin{equation}
  \left\{
  \begin{aligned}
   C_1 \cup C_3 &=\{(X_1,Y_1,Z_1)|X_1=Y_1=0\} \cup \{(X_3,Y_3,Z_3)|Y_3=Z_3=0\} \cong {\Bbb P}^1,\\
   P_2 &=\{(X_2,Y_2,Z_2)|X_2=Y_2=Z_2=0\},\\
   P_4 &=\{(X_1,Y_1,Z_1)|X_1=Z_1=0, \ Y_1=1\},\\
   P_5 &=\{(X_1,Y_1,Z_1)|X_1=0, \ Y_1=3, \ Z_1=-1\},
   \end{aligned}
  \right. 
\end{equation}
where $(X_i,Y_i,Z_i)$ are given by \eqref{P2cover}.
\end{lemma}
This lemma can be proven by a direct calculation. \qed

Next let us calculate its local index at each point.
\begin{center}
\begin{tabular}{|c|c|c|} \hline 
Singular point & Type of local index   \\ \hline 
$P_1$ & $(-1,-3,0)$  \\ \hline 
$P_2$ & $(-3,-1,-1)$  \\ \hline 
$P_4$ & $(1,3,2)$  \\ \hline 
$P_5$ & $(-1,3,-6)$  \\ \hline 
$P_6$ & $(0,-3c,-c)$  \\ \hline 
\end{tabular}
\end{center}
Here, the notations $P_1,P_3$ and $P_6$ are given by
\begin{equation}
  \left\{
  \begin{aligned}
   P_1 &=\{(X_1,Y_1,Z_1)|X_1=Y_1=Z_1=0\} \in C_1,\\
   P_3 &=\{(X_3,Y_3,Z_3)|X_3=Y_3=Z_3=0\} \in C_3,\\
   P_6 &=\{(X_3,Y_3,Z_3)|Y_3=Z_3=0, \ X_3=c\} \in C_3.
   \end{aligned}
  \right. 
\end{equation}

\begin{theorem}
The phase space ${\mathcal X}$ for the system \eqref{system:Chazy-XI} is obtained by gluing six copies of ${\Bbb C}^3${\rm:\rm}
\begin{center}
${U_j}={\Bbb C}^3 \ni \{(x_j,y_j,z_j)\},  \ \ j=0,1,\ldots,5$
\end{center}
via the following birational transformations{\rm:\rm}
\begin{align}
\begin{split}
0) \ &x_0=x, \quad y_0=y, \quad z_0=z,\\
1) \ &x_1=\frac{1}{x}, \quad y_1=x^2y, \quad z_1=\frac{z}{x},\\
2) \ &x_2=(x+z)y^2, \quad y_2=\frac{1}{y}, \quad z_2=z,\\
3) \ &x_3=\frac{x}{z}, \quad y_3=yz^2, \quad z_3=\frac{1}{z},\\
4) \ &x_4=\frac{1}{x}, \quad y_4=(y-x+2z)x^2, \quad z_4=xz,\\
5) \ &x_5=\frac{1}{x}, \quad y_5=\frac{1}{x^2y}, \quad z_5=(x+z)x^3y^2.
\end{split}
\end{align}
\end{theorem}

\section{Appendix F}

In 1881, Halphen studied an integrable third-order ordinary differential system \cite{26,27}:
\begin{equation}\label{system;C1}
  \left\{
  \begin{aligned}
   \frac{dx}{dt} &=x^2+\gamma(x-y)^2+\beta(z-x)^2+\alpha(y-z)^2,\\
   \frac{dy}{dt} &=y^2+\gamma(x-y)^2+\beta(z-x)^2+\alpha(y-z)^2,\\
   \frac{dz}{dt} &=z^2+\gamma(x-y)^2+\beta(z-x)^2+\alpha(y-z)^2,
   \end{aligned}
  \right. 
\end{equation}
where $x,y,z$ denote unknown complex variables and $\alpha,\beta,\gamma$ are complex parameters.

It is known that this system can be solved by hypergeometric functions (see \cite{26,27}).

\begin{proposition}
This system is invariant under the following transformations{\rm:\rm}
\begin{align}
\begin{split}
s_0(x,y,z;\alpha,\beta,\gamma) \rightarrow &(y,x,z;\beta,\alpha,\gamma),\\
s_1(x,y,z;\alpha,\beta,\gamma) \rightarrow &(z,y,x;\gamma,\beta,\alpha),\\
s_2(x,y,z;\alpha,\beta,\gamma) \rightarrow &(x,z,y;\alpha,\gamma,\beta),\\
\pi(x,y,z;\alpha,\beta,\gamma) \rightarrow &(y,z,x;\beta,\gamma,\alpha).
\end{split}
\end{align}
\end{proposition}

In \cite{25}, it is shown that when the three parameters $\alpha,\beta,\gamma$ are equal or when two of the parameters are $1/3$ this system reduced to the generalized Chazy equation which is a classically known third-order scalar polynomial ordinary differential equation:
\begin{align}
\begin{split}
XII:&\frac{d^3u}{dt^3}=2u\frac{d^2u}{dt^2}-3\left(\frac{du}{dt}\right)^2-\frac{4}{N^2-36}\left(6\frac{du}{dt}-u^2\right)^2,
\end{split}
\end{align}
where $N$ is a positive integer not equal to 1 or 6.

The general solution of the system \eqref{system;1} is densely branched for generic $\alpha,\beta,\gamma$ and so does not pass the Painlev\'e property.

In this section, we study the system \eqref{system;C1} from the viewpoint of its accessible singularities and local index.

The following Lemma shows that this rational vector field $\tilde v$ associated with the system \eqref{system;C1} has seven accessible singular points on the boundary divisor ${\mathcal H} \subset {\Bbb P}^3$.
\begin{lemma}
The rational vector field $\tilde v$ associated with the system \eqref{system;C1} has seven accessible singular points{\rm : \rm}
\begin{equation}
  \left\{
  \begin{aligned}
   P_1 =&\{(X_1,Y_1,Z_1)|X_1=0, \ Y_1=1, \ Z_1=1\},\\
   P_2 =&\left\{(X_1,Y_1,Z_1)|X_1=0, \ Y_1=1, \ Z_1=\frac{2\alpha+2\beta+1-\sqrt{4\alpha+4\beta+1}}{2(\alpha+\beta)} \right\},\\
   P_3 =&\left\{(X_1,Y_1,Z_1)|X_1=0, \ Y_1=1, \ Z_1=\frac{2\alpha+2\beta+1+\sqrt{4\alpha+4\beta+1}}{2(\alpha+\beta)} \right\},\\
   P_4 =&\left\{(X_1,Y_1,Z_1)|X_1=0, \ Y_1=\frac{2\alpha+2\gamma+1-\sqrt{4\alpha+4\gamma+1}}{2(\alpha+\gamma)}, \ Z_1=1 \right\},\\
   P_5 =&\left\{(X_1,Y_1,Z_1)|X_1=0, \ Y_1=\frac{2\alpha+2\gamma+1+\sqrt{4\alpha+4\gamma+1}}{2(\alpha+\gamma)}, \ Z_1=1 \right\},\\
   P_6 =&\{(X_1,Y_1,Z_1)|X_1=0, \ Y_1=\frac{2\beta+2\gamma+1-\sqrt{4\beta+4\gamma+1}}{2(\beta+\gamma)},\\
   &Z_1=\frac{2\beta+2\gamma+1-\sqrt{4\beta+4\gamma+1}}{2(\beta+\gamma)}\},\\
   P_7 =&\{(X_1,Y_1,Z_1)|X_1=0, \ Y_1=\frac{2\beta+2\gamma+1+\sqrt{4\beta+4\gamma+1}}{2(\beta+\gamma)},\\
   &Z_1=\frac{2\beta+2\gamma+1+\sqrt{4\beta+4\gamma+1}}{2(\beta+\gamma)}\},
   \end{aligned}
  \right. 
\end{equation}
where $(X_i,Y_i,Z_i)$ are given by \eqref{P2cover}.
\end{lemma}
This lemma can be proven by a direct calculation. \qed

Next let us calculate its local index at each point.
\begin{center}
\begin{tabular}{|c|c|c|} \hline 
Singular point & Local index $(a_1^{(i)},a_2^{(i)},a_3^{(i)})$   \\ \hline 
$P_1$ & $(-1,1,1)$  \\ \hline 
$P_2$ & $(-\frac{4\alpha+4\beta+1-\sqrt{4\alpha+4\beta+1}}{2(\alpha+\beta)},\frac{2}{1+\sqrt{4\alpha+4\beta+1}},-\frac{4\alpha+4\beta+1-\sqrt{4\alpha+4\beta+1}}{2(\alpha+\beta)})$  \\ \hline 
$P_3$ & $(-\frac{4\alpha+4\beta+1+\sqrt{4\alpha+4\beta+1}}{2(\alpha+\beta)},-\frac{2}{-1+\sqrt{4\alpha+4\beta+1}},-\frac{4\alpha+4\beta+1+\sqrt{4\alpha+4\beta+1}}{2(\alpha+\beta)})$  \\ \hline 
$P_4$ & $(-\frac{4\alpha+4\gamma+1-\sqrt{4\alpha+4\gamma+1}}{2(\alpha+\gamma)},-\frac{4\alpha+4\gamma+1-\sqrt{4\alpha+4\gamma+1}}{2(\alpha+\gamma)},\frac{2}{1+\sqrt{4\alpha+4\gamma+1}})$  \\ \hline 
$P_5$ & $(-\frac{4\alpha+4\gamma+1+\sqrt{4\alpha+4\gamma+1}}{2(\alpha+\gamma)},-\frac{4\alpha+4\gamma+1+\sqrt{4\alpha+4\gamma+1}}{2(\alpha+\gamma)},-\frac{2}{-1+\sqrt{4\alpha+4\gamma+1}})$  \\ \hline 
$P_6$ & $(-\frac{4\beta+4\gamma+1-\sqrt{4\beta+4\gamma+1}}{2(\beta+\gamma)},-\frac{4\beta+4\gamma+1-\sqrt{4\beta+4\gamma+1}}{2(\beta+\gamma)},-\frac{2}{1+\sqrt{4\beta+4\gamma+1}})$  \\ \hline 
$P_7$ & $(-\frac{4\beta+4\gamma+1+\sqrt{4\beta+4\gamma+1}}{2(\beta+\gamma)},\frac{1+\sqrt{4\beta+4\gamma+1}}{2(\beta+\gamma)},-\frac{4\beta+4\gamma+1+\sqrt{4\beta+4\gamma+1}}{2(\beta+\gamma)})$  \\ \hline 
\end{tabular}
\end{center}

It is easy to see that the system \eqref{system;C1} admits a rational solution:
\begin{equation}
x(t)=-\frac{1}{t-t_0}, \quad y(t)=-\frac{1}{t-t_0}, \quad z(t)=-\frac{1}{t-t_0} \quad (t_0 \in {\Bbb C}),
\end{equation}
which passes through $P_1$.

Let us take the coordinate system $(p,q,r)$ centered at the point $P_6$:
\begin{align}
\begin{split}
&p=\frac{1}{x}, \quad q=\frac{y}{x}-\frac{2\beta+2\gamma+1-\sqrt{4\beta+4\gamma+1}}{2(\beta+\gamma)},\\
&r=\frac{z}{x}-\frac{2\beta+2\gamma+1-\sqrt{4\beta+4\gamma+1}}{2(\beta+\gamma)}.
\end{split}
\end{align}
Making a linear transformation to arrive at
\begin{equation*}
\begin{pmatrix}
             \frac{dX}{dt} \\
             \frac{dY}{dt} \\
             \frac{dZ}{dt}
             \end{pmatrix}=\frac{1}{X}\left\{\begin{pmatrix}
             -\frac{4\beta+4\gamma+1-\sqrt{4\beta+4\gamma+1}}{2(\beta+\gamma)} & 0 & 0 \\
             0 & -\frac{4\beta+4\gamma+1-\sqrt{4\beta+4\gamma+1}}{2(\beta+\gamma)} & K_1 \\
             0 & 0 & -\frac{2}{1+\sqrt{4\beta+4\gamma+1}}
             \end{pmatrix}\begin{pmatrix}
             X \\
             Y \\
             Z 
             \end{pmatrix}+\cdots\right\},
             \end{equation*}
where $K_1$ is given by
\begin{align}
\begin{split}
K_1:=&-\frac{2\beta+2\gamma+1-\sqrt{4\beta+4\gamma+1}}{2(\beta+\gamma)}\\
&-\frac{\sqrt{2}(\beta+\gamma)}{\sqrt{1+4(\beta+\gamma)+2(\beta+\gamma)^2+(2\beta+2\gamma+1)\sqrt{4\beta+4\gamma+1}}}.
\end{split}
\end{align}

Let us take the coordinate system $(p,q,r)$ centered at the point $P_7$:
\begin{align}
\begin{split}
&p=\frac{1}{x}, \quad q=\frac{y}{x}-\frac{2\beta+2\gamma+1+\sqrt{4\beta+4\gamma+1}}{2(\beta+\gamma)},\\
&r=\frac{z}{x}-\frac{2\beta+2\gamma+1+\sqrt{4\beta+4\gamma+1}}{2(\beta+\gamma)}.
\end{split}
\end{align}
Making a linear transformation to arrive at
\begin{equation*}
\begin{pmatrix}
             \frac{dX}{dt} \\
             \frac{dY}{dt} \\
             \frac{dZ}{dt} 
             \end{pmatrix}=\frac{1}{X}\left\{\begin{pmatrix}
             -\frac{4\beta+4\gamma+1+\sqrt{4\beta+4\gamma+1}}{2(\beta+\gamma)} & 0 & 0 \\
             0 & \frac{1+\sqrt{4\beta+4\gamma+1}}{2(\beta+\gamma)} & 0 \\
             0 & K_2 & -\frac{4\beta+4\gamma+1+\sqrt{4\beta+4\gamma+1}}{2(\beta+\gamma)}
             \end{pmatrix}\begin{pmatrix}
             X \\
             Y \\
             Z 
             \end{pmatrix}+\cdots\right\},
             \end{equation*}
where $K_2$ is given by
\begin{align}
\begin{split}
K_2:=&-\frac{2\beta+2\gamma+1+\sqrt{4\beta+4\gamma+1}}{2(\beta+\gamma)}\\
&-\frac{\sqrt{2}(\beta+\gamma)}{\sqrt{1+4(\beta+\gamma)+2(\beta+\gamma)^2-(2\beta+2\gamma+1)\sqrt{4\beta+4\gamma+1}}}.
\end{split}
\end{align}

\begin{proposition}
Each local index at each accessible singular point $P_i, \ (i=1,2,\ldots,7)$ satisfies the condition$:$
\begin{equation}
\left(\frac{a_2^{(i)}}{a_1^{(i)}},\frac{a_3^{(i)}}{a_1^{(i)}} \right) \in {\Bbb Z}^2
\end{equation}
if and only if the parameters $\alpha,\beta$ and $\gamma$ satisfy the conditions$:$
\begin{equation}
  \left\{
  \begin{aligned}
   \frac{1}{\sqrt{4\alpha+4\beta+1}} &=l,\\
   \frac{1}{\sqrt{4\alpha+4\gamma+1}} &=m,\\
   \frac{1}{\sqrt{4\beta+4\gamma+1}} &=n,
   \end{aligned}
  \right. 
\end{equation}
where $(l,m,n) \in {\Bbb Z}^3$. 
\end{proposition}
This equation can be solved by
\begin{equation}\label{aaaaa}
  \left\{
  \begin{aligned}
   \alpha &=\frac{1}{8}\left(\frac{1}{l^2}+\frac{1}{m^2}-\frac{1}{n^2}-1 \right),\\
   \beta &=\frac{1}{8}\left(\frac{1}{l^2}-\frac{1}{m^2}+\frac{1}{n^2}-1 \right),\\
   \gamma &=\frac{1}{8}\left(-\frac{1}{l^2}+\frac{1}{m^2}+\frac{1}{n^2}-1 \right),
   \end{aligned}
  \right. 
\end{equation}
where $(l,m,n) \in {\Bbb N}^3$. Under the condition \eqref{aaaaa}, we see that
\begin{equation*}
K_1=0, \quad K_2=0.
\end{equation*}

Next, we give a generalization in dimension four of the second Halphen equation given by
\begin{equation}\label{system;a1}
  \left\{
  \begin{aligned}
   \frac{dx}{dt} &=x^2+\alpha(x-y)^2+\beta(x-z)^2+\chi(x-w)^2+\delta(y-z)^2+\varepsilon(y-w)^2+\gamma(z-w)^2,\\
   \frac{dy}{dt} &=y^2+\alpha(x-y)^2+\beta(x-z)^2+\chi(x-w)^2+\delta(y-z)^2+\varepsilon(y-w)^2+\gamma(z-w)^2,\\
   \frac{dz}{dt} &=z^2+\alpha(x-y)^2+\beta(x-z)^2+\chi(x-w)^2+\delta(y-z)^2+\varepsilon(y-w)^2+\gamma(z-w)^2,\\
   \frac{dw}{dt} &=w^2+\alpha(x-y)^2+\beta(x-z)^2+\chi(x-w)^2+\delta(y-z)^2+\varepsilon(y-w)^2+\gamma(z-w)^2,\\
   \end{aligned}
  \right. 
\end{equation}
where $x,y,z,w$ denote unknown complex variables and $\alpha,\beta,\chi,\delta,\varepsilon,\gamma$ are complex parameters.

\begin{proposition}
This system is invariant under the following transformations{\rm:\rm}
\begin{align}
\begin{split}
s_1(*) \rightarrow &(y,x,z,w;\alpha,\delta,\varepsilon,\beta,\chi,\gamma),\\
s_2(*) \rightarrow &(z,y,x,w;\delta,\beta,\gamma,\alpha,\varepsilon,\chi),\\
s_3(*) \rightarrow &(w,y,z,x;\varepsilon,\gamma,\chi,\delta,\alpha,\beta),\\
s_4(*) \rightarrow &(x,z,y,w;\beta,\alpha,\chi,\delta,\gamma,\varepsilon),\\
s_5(*) \rightarrow &(x,w,z,y;\chi,\beta,\alpha,\gamma,\varepsilon,\delta),\\
s_6(*) \rightarrow &(x,y,w,z;\alpha,\chi,\beta,\varepsilon,\delta,\gamma),\\
\pi(*) \rightarrow &(y,z,w,x;\delta,\varepsilon,\alpha,\gamma,\beta,\chi),
\end{split}
\end{align}
where the symbol $(*)$ denotes $(x,y,z,w;\alpha,\beta,\chi,\delta,\varepsilon,\gamma)$, and $s_i^2=1, \ {\pi}^4=1$.
\end{proposition}

\begin{lemma}
The system \eqref{system;a1} has fifteen accessible singular points on the boundary divisor ${\mathcal H} \subset {\Bbb P}^4$.
\end{lemma}

It is easy to see that the system \eqref{system;a1} admits a rational solution:
\begin{equation}
x(t)=-\frac{1}{t-t_0}, \quad y(t)=-\frac{1}{t-t_0}, \quad z(t)=-\frac{1}{t-t_0}, \quad w(t)=-\frac{1}{t-t_0} \quad (t_0 \in {\Bbb C}).
\end{equation}

For the system \eqref{system;a1}, each local index at each accessible singular point $P_i, \ (i=1,2,\ldots,15)$ satisfies the condition:
\begin{equation}
\left(\frac{a_2^{(i)}}{a_1^{(i)}},\frac{a_3^{(i)}}{a_1^{(i)}},\frac{a_4^{(i)}}{a_1^{(i)}} \right) \in {\Bbb Z}^3
\end{equation}
if and only if the parameters $\alpha,\beta,\chi,\delta,\varepsilon$ and $\gamma$ satisfy the conditions:
\begin{equation}
  \left\{
  \begin{aligned}
   4\beta+4\delta+4\gamma+1 &=\frac{1}{m_1^2},\\
   4\chi+4\varepsilon+4\gamma+1 &=\frac{1}{m_2^2},\\
   4\alpha+4\delta+4\varepsilon+1 &=\frac{1}{m_3^2},\\
   4\beta+4\chi+4\delta+4\varepsilon+1 &=\frac{1}{m_4^2},\\
   4\alpha+4\beta+4\chi+1 &=\frac{1}{m_5^2},\\
   4\alpha+4\beta+4\varepsilon+4\gamma+1 &=\frac{1}{m_6^2},\\
   4\alpha+4\chi+4\delta+4\gamma+1 &=\frac{1}{m_7^2},
   \end{aligned}
  \right. 
\end{equation}
where $m_i \in {\Bbb Z} \ (i=1,2,\ldots,7)$. This equation can be solved by
\begin{equation}
  \left\{
  \begin{aligned}
   \alpha &=\frac{1}{8}\left(\frac{1}{m_3^2}-\frac{1}{m_4^2}+\frac{1}{m_5^2}-1 \right),\\
   \beta &=\frac{1}{8}\left(\frac{1}{m_1^2}+\frac{1}{m_5^2}-\frac{1}{m_7^2}-1 \right),\\
   \chi &=\frac{1}{8}\left(-\frac{1}{m_1^2}-\frac{1}{m_3^2}+\frac{1}{m_4^2}+\frac{1}{m_7^2} \right),\\
   \delta &=\frac{1}{8}\left(-\frac{1}{m_2^2}+\frac{1}{m_4^2}-\frac{1}{m_5^2}+\frac{1}{m_7^2} \right),\\
   \varepsilon &=\frac{1}{8}\left(\frac{1}{m_1^2}+\frac{1}{m_2^2}-\frac{1}{m_4^2}-1 \right),\\
   \gamma &=\frac{1}{8}\left(\frac{1}{m_2^2}+\frac{1}{m_3^2}-\frac{1}{m_7^2}-1 \right)\\
   \end{aligned}
  \right. 
\end{equation}
with the condition
\begin{equation}
-\frac{1}{m_7^2}-\frac{1}{m_6^2}+\frac{1}{m_5^2}-\frac{1}{m_4^2}+\frac{1}{m_3^2}+\frac{1}{m_2^2}+\frac{1}{m_1^2}=1 \quad (m_i \in {\Bbb N}).
\end{equation}

It is still an open question whether the system \eqref{system;a1} can be solved by known functions.

\end{document}